\documentclass[a4paper,7pt,oneside,onecolumn,number,preprint,centertitle]{elsarticle}

\usepackage{amsmath,amssymb,bm}
\usepackage{graphicx}
\usepackage{textcomp}
\usepackage{xcolor}
\usepackage{amsthm}
\usepackage{tikz}

\usepackage{enumerate}
\usepackage{bbold}
\usepackage{epsfig}
\usepackage{multicol}
\usepackage{float}
\usepackage{hyperref}
\usepackage{multirow}
\usepackage{mathrsfs} %

\usepackage{algorithm}
\usepackage{algorithmicx}
\usepackage{algpseudocode}
\usepackage{flushend}

\usepackage{ifthen}
\DeclareMathAlphabet{\mbb}{U}{bbold}{m}{n} %
\newcommand \mb[1] {\ifthenelse{\equal{#1}{0}}{\mbb{0}}{\ifthenelse{\equal{#1}{1}}{\mbb{1}}{\mathbb{#1}}}} %

\newcommand \mc {\mathcal}

\newcommand \RR {\mb R}

\newcommand \T {^\top}

\renewcommand \d {\mathrm{d}}

\newcommand \diag {\mathrm{diag}}

\newcommand{\on}[1]{\operatorname{#1}}

\renewcommand \subset {\subseteq}

\renewcommand{\l}{\left}

\newtheorem{thm}{Theorem}
\newtheorem{prop}[thm]{Proposition}
\newtheorem{cor}[thm]{Corollary}
\newtheorem{rem}{Remark}
\newtheorem{lem}[thm]{Lemma}
\newtheorem{ass}{Assumption}

\newtheorem{exmp}{Example}%

\renewcommand{\textcolor}[2]{#2}
\begin{document}
\begin{frontmatter}

\title{Minimax Linear Regulator Problems for \\ Positive Systems}

\author[inst1]{Alba Gurpegui}
    \ead{alba.gurpegui_ramon@control.lth.se}

\affiliation[inst1]{organization={Department of Automatic Control, Lund University},%
            addressline={Ole Römers väg 1}, 
            city={Lund},
            postcode={223 63}, 
            country={Sweden}}

\author[inst1]{Mark Jeeninga}
    \ead{mark.jeeninga@control.lth.se}
    
\author[inst1]{Emma Tegling}
    \ead{emma.tegling@control.lth.se}
\author[inst1]{Anders Rantzer}
    \ead{anders.rantzer@control.lth.se}

\begin{abstract}
Explicit solutions to optimal control problems are rarely obtainable. Of particular interest are the explicit solutions derived for minimax problems, providing a framework to address adversarial conditions and uncertainty. This work considers a multi-disturbance minimax Linear Regulator (LR) framework for positive linear time-invariant systems in continuous time, which, analogous to the Linear-Quadratic Regulator (LQR) problem, can be utilized for the stabilization of positive systems.
The problem is studied for nonnegative and state-bounded disturbances.
Dynamic programming theory is leveraged to derive explicit solutions to the minimax LR problem for both finite and infinite time horizons.
In addition, a fixed-point method is proposed that computes the solution for the infinite horizon case, and the minimum $L_1$-induced gain of the system is studied.

We motivate the prospective scalability properties of our framework with a large-scale water management network.

\footnotetext{
The authors are members of the ELLIIT Strategic Research Area at Lund University. This work is partially funded by Wallenberg AI, Autonomous Systems and Software Program (WASP) and the European Research Council (ERC) under the European Union's Horizon 2020 research and innovation programme under grant agreement No 834142  (ScalableControl).}
\end{abstract}


\end{frontmatter}

\section{Introduction}
\label{sec:introduction}
\subsection{Motivation}\label{subsec:motiv}
Minimax optimal control problems are widespread across control theory and engineering disciplines. They offer a methodology for formulating and addressing challenges characterized by competitive elements and uncertainties. Their formulation often involves solving the Hamilton-Jacobi-Isaacs equation~\cite{HJI}, a fundamental partial differential equation
governing differential games and worst-case control strategies. 
Such problems are prevalent in areas including robust control, multi-agent systems and game theory, where it finds direct application in $H_{\infty}$ robust control~\cite{aux16, aux17, aux18, aux19, aux20, aux21}. Tackling these problems presents significant difficulties, particularly when we deal with large-scale systems due to the computational complexity and dimensionality of the associated optimization problems. A generic nonnegative linear cost minimax setting is studied to design control systems that are robust to worst-case uncertainties and disturbances.

The study of optimal control problems with nonnegative cost dates back to~\cite{linear_cost1, linear_cost2, linear_cost3}. %
In particular, this work draws inspiration from the unified dynamic programming framework developed by D. Bertsekas \cite{Bertsekas_dynprograndoptcontrol, BertsekasVIPI}. Despite significant advancements in the theory of dynamic programming, obtaining exact solutions often remains challenging. A notable exception is the linear quadratic problem, as presented in the pioneering work of Kalman \cite{Kalman}. In both quadratic and linear nonnegative cost settings in infinite horizon, the optimal cost is determined by an algebraic equation. For large-scale systems, when we choose the linear quadratic approach, the algebraic equation is the Riccati equation, and the number of unknown parameters grows quadratically with the state dimension. In the linear regulator setting, the algebraic equation derived in ~\cite{yuchao} for the optimal cost is instead linear with respect to the dimension of the state. %
This difference becomes 
significant when the state dimension becomes large. 

Another important characteristic of our problem setup is the positive dynamics. Positive systems are characterized by the property that their state and output remain nonnegative for any nonnegative input and initial state. This type of dynamics may capture a wide range of physical phenomena, and they have thus been the subject of substantial control theoretic research. Classical references on this topic involve the well-established positive matrix theory pioneered by Perron~\cite{Perron1907} and Frobenius~\cite{Frobenius}, summarized in classical textbooks~\cite{BermanBook} and~\cite{Luenberger}. A primary advantage of positive systems is that stability can be easily verified using linear Lyapunov functions~\cite{Blanchini_lyapunov, blanchini_lyapunov2}. This is especially important in the large scale setting because of the resulting computational scalability~\cite{EbiharaPeucelleArzelierTAC2017,tutorial}.  The study of optimal control over finite time horizons naturally leads to switching control laws, creating an organic link to the rich literature on switched positive systems~\cite{aux9, aux10, 12tutorial}. In particular, the worst-case disturbance for a bilinear system with an external positive disturbance has been characterized in recent work on compartmental models~\cite{BLANCHINI_REV2}. The present paper has several connections to that work, which will be revisited throughout the manuscript. Additionally, the $L_{1}$-induced gain of a positive system plays an important role in robust stability analysis against dynamical and parametric uncertainties~\cite{Ebihara_robustnessl1, tutorial, Briat2013, Ebihara1}. This has led to contributions dealing with
decentralized and distributed control for large-scale positive systems~\cite{ebihara_extra1,ebihara_extra2}.

Inspired by the novel class of optimal control problems involving nonnegative linear costs~\cite{yuchao}, previous research in discrete time~\cite{AlbaEmmaAnders} presents an explicit solution for a class of minimax optimal control problems with linear objective functions, positive linear systems, and homogeneous constraints on control and disturbances. 
In the current paper, these results are extended to a minimax setting in continuous time with multiple disturbances. Two types of disturbances are considered: bounded by elementwise linear constraints, and unconstrained nonnegative disturbances. 

The main contributions of this manuscript are as follows: 

\begin{enumerate}[C1.]
    \item Necessary conditions on the parameters of the minimax control problem are formulated, for which, if the optimal solution is finite, the closed-loop system is a positive system (Theorems~\ref{theorem_finite},~\ref{MainThm_minimax_cont_inf}) and the cost is minimized.
    \item Dynamic programming is applied without imposing predefined constraints on linearity or sparsity; instead, these properties arise naturally from the optimization criteria and constraints. This approach differs from methods in which the control design explicitly dictates specific structures.
    \item In both finite and infinite horizon settings, the minimization and the two distinct maximization problems are decoupled. This property enables the explicit solution of the HJI equation to be extended to the minimax framework. The HJI equation takes the form of an ordinary differential equation (ODE) in the finite horizon setting, whereas it reduces to an algebraic equation in the infinite horizon case.
    \item A fixed-point method for computing the solution to the algebraic HJI equation in the presence of elementwise bounded disturbances (Theorem~\ref{Cont_VI}).
    \item 
    The stabilizability and detectability of the Linear Regulator (LR) problem are studied, and a priori detectability conditions of the systems under consideration are established (Section~\ref{Stab_Detect_subsec}).
    \item 
     A linear programming formulation is proposed to solve the Isaacs equation of the LR problem under appropriate assumptions. Moreover, an analysis of the primal and dual linear program formulations are provided and necessary conditions are presented for which the LR problem stabilizes the system (Section~\ref{LP_subsect}). 
    \item Analysis of the $L_1$-induced gain of the system and its relationship to the disturbance penalty in the cost function of the minimax setting is given (Section~\ref{section_l1gain}).
\end{enumerate}
The present manuscript consolidates and extends prior results in discrete time~\cite{yuchao, AlbaEmmaAnders}, into a unified continuous-time framework for minimax linear regulator problems. Specifically, Contributions C1–C3 establish the continuous-time dynamic programming formulation and explicit solutions to the associated Hamilton–Jacobi–Isaacs equations, thereby generalizing the discrete-time minimax framework of~\cite{AlbaEmmaAnders} to a full dynamic game setting with multiple disturbances. Contributions C4–C6 further extend the discrete-time linear regulator theory of~\cite{yuchao} by providing its continuous-time counterpart, complete with a fixed-point computational approach and linear programming characterization. Together, these results form a continuous-time framework that encompasses robustness, positivity, and scalability properties previously developed in discrete-time formulations. Finally, C7 offers a new interpretation of the $L_1$-induced gain within this framework, linking disturbance attenuation directly to the minimax optimization structure. We note that a related problem concerning maximal disturbances has been studied in \cite{BLANCHINI_REV2}.
\subsection{Problem Setup}\label{subsec_setup}
Consider the LTI system
\begin{align}\label{cont_syst}
    \dot{x}(t)=Ax(t)+ Bu(t)+Fw(t) +Hv(t)
\end{align}
where $x$ represents the $n-$dimensional vector of state variables, $u$ the $m-$dimensional control variable, $w$ the $l-$dimensional disturbance, $v$ the $c-$dimensional disturbance, $A \in \mathbb{R}^{n\times n}$ is Metzler, $B \in \mathbb{R}^{n\times m }$, $F \in \mathbb{R}^{n \times l}$ and $H \in \mathbb{R}^{n \times c}$. In this work, the control $u$ and the disturbance $v$ are assumed to be bounded variables, whereas the disturbance $w$ is only assumed to be nonnegative.

The general state-space model considered in this paper is,
\begin{align}\label{dynamics}
   G_{\mu}: \left\{\begin{aligned}
\dot{x}(t)&=Ax(t)+B\mu(x(t))+Fw(t)+Hv(t)\\ 
z(t)&=s^{\top}x(t)+r^{\top}\mu(x(t)) - \boldsymbol{\gamma}^{\top} w-\delta^{\top}v 
\end{aligned}\right.
\end{align}
 where $\mu$ is any, potentially nonlinear, control policy, $z$ is a target output variable representing the system's performance, $s \in \mathbb{R}^{n} $ and $r \in \mathbb{R}^{m} $. Assuming zero initial conditions, the system dynamics~\eqref{dynamics} can be seen as an operator $G_{\mu}$ from the disturbance $w$ to the output $z$.

 The objective is to minimize the worst-case cost over all possible control strategies. The cost is the integral of the system's output $z$, which is a linear function. %
 Linear cost functions are suitable for positive systems, where the variables that may represent e.g. flow, population, or inventory, are inherently positive. It is reasonable to consider a linear performance output, as this directly penalizes deviations, which naturally aligns with the systems' physical constraints.
 This approach also provides a clear interpretation of costs related to resource usage or penalties for exceeding limits.
\begin{figure}[t!]
\centering
\tikzset{every picture/.style={line width=0.75pt}} %

\begin{tikzpicture}[x=0.75pt,y=0.75pt,yscale=-1,xscale=1]
\draw   (110.12,74.85) -- (279.51,74.85) -- (279.51,123) -- (110.12,123) -- cycle ;
\draw    (366.5,85.61) -- (282.09,85.61) ;
\draw [shift={(280.09,85.61)}, rotate = 360] [color={rgb, 255:red, 0; green, 0; blue, 0 }  ][line width=0.75]    (10.93,-3.29) .. controls (6.95,-1.4) and (3.31,-0.3) .. (0,0) .. controls (3.31,0.3) and (6.95,1.4) .. (10.93,3.29)   ;
\draw    (110.69,98.54) -- (32,98.54) ;
\draw [shift={(30,98.54)}, rotate = 360] [color={rgb, 255:red, 0; green, 0; blue, 0 }  ][line width=0.75]    (10.93,-3.29) .. controls (6.95,-1.4) and (3.31,-0.3) .. (0,0) .. controls (3.31,0.3) and (6.95,1.4) .. (10.93,3.29)   ;
\draw   (37.82,215.49) .. controls (37.81,220.16) and (40.13,222.5) .. (44.8,222.51) -- (188.61,222.97) .. controls (195.28,222.99) and (198.6,225.33) .. (198.59,230) .. controls (198.6,225.33) and (201.94,223.01) .. (208.61,223.03)(205.61,223.02) -- (352.42,223.48) .. controls (357.09,223.49) and (359.43,221.17) .. (359.44,216.5) ;
\draw   (170.21,170.14) -- (215.42,170.14) -- (215.42,203.26) -- (170.21,203.26) -- cycle ;
\draw    (326.61,116.55) -- (283.85,116.66) ;
\draw [shift={(281.85,116.66)}, rotate = 359.86] [color={rgb, 255:red, 0; green, 0; blue, 0 }  ][line width=0.75]    (10.93,-3.29) .. controls (6.95,-1.4) and (3.31,-0.3) .. (0,0) .. controls (3.31,0.3) and (6.95,1.4) .. (10.93,3.29)   ;
\draw    (66.63,189.5) -- (167.1,188.73) ;
\draw [shift={(169.1,188.71)}, rotate = 179.56] [color={rgb, 255:red, 0; green, 0; blue, 0 }  ][line width=0.75]    (10.93,-3.29) .. controls (6.95,-1.4) and (3.31,-0.3) .. (0,0) .. controls (3.31,0.3) and (6.95,1.4) .. (10.93,3.29)   ;
\draw    (66.63,115.62) -- (66.63,189.5) ;
\draw    (326.61,116.55) -- (327.59,190.01) ;
\draw    (66.63,115.62) -- (108.76,115.46) ;
\draw    (215.13,189.81) -- (327.59,190.01) ;
\draw    (366.5,96.7) -- (282.09,96.7) ;
\draw [shift={(280.09,96.7)}, rotate = 360] [color={rgb, 255:red, 0; green, 0; blue, 0 }  ][line width=0.75]    (10.93,-3.29) .. controls (6.95,-1.4) and (3.31,-0.3) .. (0,0) .. controls (3.31,0.3) and (6.95,1.4) .. (10.93,3.29)   ;

\draw (107.45,79.39) node [anchor=north west][inner sep=0.95pt]  [font=\small]  {$ \begin{array}{l}
\hspace{2mm}\dot{x} =Ax+Bu+Fw+Hv\\
z=s^{\top } x+r^{\top } u-\boldsymbol\gamma^{\top}w-\delta^{\top}v
\end{array}$};
\draw (337.95,70.58) node [anchor=north west][inner sep=0.75pt]  [font=\normalsize]  {$w$};
\draw (69.21,79.58) node [anchor=north west][inner sep=0.75pt]  [font=\normalsize]  {$z$};
\draw (182.81,234.46) node [anchor=north west][inner sep=0.75pt]  [font=\large]  {$G_{\mu }$};
\draw (182.76,176.68) node [anchor=north west][inner sep=0.75pt]  [font=\large]  {$\mu ( \cdot )$};
\draw (338.97,101.88) node [anchor=north west][inner sep=0.75pt]  [font=\normalsize]  {$v$};
\draw (339.93,141.29) node [anchor=north west][inner sep=0.75pt]  [font=\normalsize]  {$u$};
\draw (46.39,143.48) node [anchor=north west][inner sep=0.75pt]  [font=\normalsize]  {$x$};
\end{tikzpicture}

\normalsize\caption{Block Diagram of the closed-loop system dynamics in~\eqref{dynamics} under the presence of two types of disturbances.}
\label{diagram}
\end{figure}
\subsection{Notation and Preliminaries}\label{subsect_notation} Let $\mathbb{R}_{+}$ denote the set of nonnegative real numbers, $\mathbb{R}^{n}$ the $n-$dimensional Euclidean space and $\mathbb{R}^{n}_{+}$ the positive orthant. Let $\mathbb{R}^{m \times n}$ denote the set of $m$ by $n$ matrices. Any vector is, by default, a column vector. Row vectors are specified explicitly. We use $\mathbb{0}$ to denote the vector/matrix with all zero elements and $\mathbb{1}$ to denote a column vector with unit entries in $\mathbb{R}$ of appropriate dimensions. For a real matrix $X$ we use the notation $\left | X \right |$ to denote the matrix obtained by replacing the elements of $X$ with their absolute values. A real matrix $X$ is called nonnegative $X \geq 0$ (respectively nonpositive $X\leq0)$ if all the elements of $X$ are nonnegative (nonpositive). If all the elements of $X$ are strictly positive (strictly negative), we call $X$ a positive $X>0$ (negative $X<0$) matrix. For real matrices $X,Y$, the inequality $X \geq Y$ $(X \leq Y)$ means that all the elements of the matrix $X-Y$ are nonnegative (nonpositive). The spectral abscissa $\alpha(X)$ of a matrix $X$ is the greatest real part of the eigenvalues of $X$. Because column and row vectors are special forms of matrices, the function $\left | \cdot  \right |$ and relations $\geq$, $\leq$, $>$, $<$ apply to them. A square matrix $A \in \mathbb{R}^{n \times n}$ is \textit{Metzler} if its off-diagonal entries are all nonnegative. A square matrix is \textit{Hurwitz} if all its eigenvalues have strictly negative real part.
The signum of a scalar is defined as the set-valued map 
\begin{align*}
    \operatorname{sign}(x) = 
    \scalebox{.9}{\ensuremath{
    \begin{cases}
        \{+1\} & \text{if } x > 0\\
        \{-1\} & \text{if } x < 0\\
        [-1,+1] & \text{if } x = 0
    \end{cases}}}
\end{align*}
and is extended to vectors in an element-wise fashion.
\section{The Minimax Linear Regulator Problem}
The optimal control problem in this section is formulated as a continuous-time minimax problem with nonnegative linear cost, positive linear dynamics, elementwise linear constraints on the control policy and the disturbance $v \in \mathbb{R}^{c}$, nonnegative unconstrained disturbance $w$ and zero terminal cost,
\begin{align}\label{minimaxProb_finite}
        &\underset{\mu}{\inf} \hspace{0.3mm}\underset{w,v}{\mathrm{max}}\hspace{1mm}\int_{0}^{T} \left [ s^{\top}x(\tau) +r^{\top}u(\tau)-\boldsymbol\gamma^{\top}w(\tau)-\delta^{\top}v(\tau)\right ]d\tau \notag\\
        &\mathrm{Subject \hspace{2mm} to} \notag \\
        &~~~~~~~~\dot{x}(t)=Ax(t)+Bu(t)+Fw(t)+Hv(t) \\
        &~~~~~~~~x(0)=x_{0},\hspace{2mm} x_0 \in \mathbb{R}^n_+ , \hspace{2mm} u(t)=\mu(x(t)), \notag\\
        &~~~~~~~~\left | u \right | \leq Ex, \hspace{2mm} w \geq 0, \hspace{2mm} \left | v \right | \leq Gx. \notag
    \end{align}
Recall the state $x \in \mathbb{R}^n$, control signal $u \in \mathbb{R}^m$, disturbances $w\in \mathbb{R}^l$ and $v\in \mathbb{R}^c$. A control policy $\mu$ is considered, which determines the control input as a function of the state. The control policy $\mu$ selects a control input $u$, which is subject to state-dependent bounds prescribed by the matrix $E$. The matrix $G$ prescribes similar bounds on the disturbance $v$. The objective is to minimize the worst case cost over all possible policies $\mu$ for the time horizon $T\in \mathbb R_+$.%

The main results of this section, Theorems~\ref{theorem_finite} and~\ref{MainThm_minimax_cont_inf}, make use of dynamic programming theory (see the Appendix~\ref{Appendix}) to derive explicit solutions to the Hamilton-Jacobi-Isaacs (HJI) equation for the finite and infinite horizon problem settings~\eqref{minimaxProb_finite} and~\eqref{minimax_prob_setup_cont_inf}, respectively. These results also provide necessary and sufficient conditions for the existence of finite solutions. We remark that~\cite{BLANCHINI_REV2} presents analogous expressions of the explicit solutions~\eqref{p_ODE} and~\eqref{p_eq_minimax_inf} for a bilinear formulation that, for some purposes, turns out to be equivalent to our problem setting (in the absence of disturbances).
\subsection{Finite Horizon}\label{subsection_finite_minimax}
The next theorem characterizes the optimal solutions to the problem setting~\eqref{minimaxProb_finite}, when the time horizon is finite.
\begin{thm}\label{theorem_finite}
    Let $A \in \mathbb{R}^{n\times n}$, $B = \left [ B_{1} \dots B_{m} \right ] \in \mathbb{R}^{n\times m }$, $F \in \mathbb{R}^{n \times l}_+$, $H \in \mathbb{R}^{n \times c}$, $E = \left [ E_{1}^{\top} \dots E_{m}^{\top} \right ]^{\top} \in \mathbb{R}_{+}^{m \times n}$ such that $E_i^{\top} \neq \mathbb{0}$ for all $i$, $G\in \mathbb{R}^{c \times n}_{+} $, $s \in \mathbb{R}^{n} $, $r \in \mathbb{R}^{m} $, $\boldsymbol\gamma \in \mathbb{R}^{l}_{+}$ and $\delta \in \mathbb{R}^{c}$. Suppose that 
    \begin{align}\label{ass_a}
    &A -\left | B \right |E -\left | H \right |G  \hspace{1.5mm} \text{Metzler},\\
    &s \geq E ^{\top} \left | r \right | - G^{\top} \left | \delta \right |.
    \label{ass_s}
    \end{align} 
    Then the following statements are equivalent:
    \begin{enumerate}[(i)]
        \item The optimal control problem~\eqref{minimaxProb_finite} has a finite optimal value for every $x_{0}\in \mathbb{R}^{n}_{+}$;
        \item The differential equation in $p(t)\in \mathbb{R}^n$,
        \begin{align}\label{p_ODE}
            &\scalebox{.9}{$-\dot{p}(t)= s + A^{\top}p(t) -E^{\top}\left | r+B^{\top}p(t) \right |+ G^{\top}\!\left | -\delta\!+ \!H^{\top}p(t) \right |,$} \notag\\
            &\hspace{1.2mm}\scalebox{.9}{$p(T)=0,$}
        \end{align}
        has a unique solution and  
 \begin{align}\label{gamma_ass_finite}
            \boldsymbol{\gamma}  \geq    F^{\top}p(0). 
\end{align}
    \end{enumerate}
Moreover, if the above conditions hold, the minimal value of the optimal control problem~\eqref{minimaxProb_finite} is $p(0)^{\top}x_{0}$ and the optimal control policy is given by $u^{*}(t)=-K(t)x^{*}(t)$ with
\begin{align} \label{K_cont_finite}
    K(t) \in \begin{bmatrix}
\mathrm{sign}(r_{1}+ p(t)^{\top}B_{1})E_{1}\\ 
\vdots \\ 
\mathrm{sign}(r_m + p(t)^{\top} B_m) E_m
\end{bmatrix}.
\end{align}
\end{thm}
We refer to the entries of the vector $r^{\top}+p(t)^{\top}B$ as the \textit{input gradients}.

\begin{rem}\label{rem_K}
    The right-hand side of \eqref{K_cont_finite} is a set, since multiple controllers may exist that achieve the optimal cost.
    For any index $i$ such that the input gradient satisfies $r_i+p(t)_i^{\top}B_i=0$ all feedback gain matrices are in the set 
\begin{align*}
    \scalebox{.92}{$\mathcal{K}=\left\{DE \hspace{1mm} | \hspace{1mm} D_{ii}\in \left[ -1, 1\right], D_{jj} \in \operatorname{sign}(r_j+p(t)_j^{\top}B_j) \hspace{1mm}\mathrm{ for } \hspace{1mm} j\neq i \right\}$}
\end{align*}    
    lead to the same and unique solution $p(t)$ of the HJI ODE equation~\eqref{p_ODE}.
\end{rem}
\begin{rem}\label{rem_ass_pos}
   The condition $A-|B|E-|H|G$ Metzler ensures invariance of the positive orthant under the system dynamics. Indeed, this assumption implies that both the open-loop and closed-loop system dynamics correspond to those of a positive system.
\end{rem}
\begin{rem}\label{rem_1}  
Condition~\eqref{ass_s} ensures that the cost function is bounded from below. This is not merely a technical consideration but also sheds light on the influence of disturbance and control penalties within the cost function. In particular, it illustrates that if the control penalty is excessively large, there is little incentive to apply control, whereas if the disturbance penalty is too low, there is insufficient motivation to counteract the disturbance.
\end{rem}

\begin{rem}\label{rem_2}
    Among all the introduced, potentially nonlinear and time-varying policies $\mu$, the optimal policy turns out to be linear. Moreover, the sparsity structure of the control gain $K(t)$ in \eqref{K_cont_finite} is inherited from the $E$ matrix which can be determined by the problem designer. 
\end{rem}

To prove Theorem~\ref{theorem_finite} we use the following Lemma.
\begin{lem}\label{lem_monotone}
    Assume that~\eqref{ass_a} and~\eqref{ass_s} hold. The solutions to the differential equation~\eqref{p_ODE} are monotone for all nonnegative terminal conditions. That is, if $p_i$, $i=1,2$ satisfy the differential equation~\eqref{p_ODE} and $0 \leq p_1(t)\leq p_2(t)$ then $p_1(\tau)\leq p_2(\tau)$ for all $\tau \le t$. Moreover, if $-\dot p(t)\geq 0$ then $-\dot p(\tau)\geq 0$ for all $\tau \leq t$.
\end{lem}
\begin{proof} 
Let $-\dot p_i(t)=f_p(p_i)$, where $f_p(p)$ is the right-hand side of \eqref{p_ODE}. By the triangle inequality it holds that $\frac{\partial f_p}{\partial p}\geq A^{\top}-E^{\top}|B^{\top}|-G^{\top}|H^{\top}|$, which is Metzler by assumption~\eqref{ass_a}, hence $\frac{\partial f_p}{\partial p}$ is Metzler. By~\cite[Rem 1.1, Ch. 3]{Smith1995MonotoneDS} the trajectories that satisfy $-\dot p = f_p(p)$ are monotone. 
Similarly, by assumptions~\eqref{ass_a},~\eqref{ass_s} and the triangle inequality it follows that equation~\eqref{p_ODE} satisfies the inequality
\begin{align*}
   -\dot{p}(t)\geq (A^{\top}-E^{\top}|B^{\top}|-G^{\top}|H^{\top}|)p(t),
\end{align*}
which, again by~\eqref{ass_a}, completes the proof.
\end{proof}

Next we prove Theorem~\ref{theorem_finite}.

\begin{proof}
 $(i)\Longrightarrow(ii)$ %
By the Picard-Lindelöf Theorem \cite{kelley2010theory} there exists a unique $p: [0,T] \rightarrow \mathbb R^n$ solving~\eqref{p_ODE}. Next, we will show that if $\boldsymbol{\gamma} \geq F^{\top}p(t)$ is violated on $I\subseteq [0,T]$, then the cost of the optimal control problem~\eqref{minimaxProb_finite} is unbounded. This shows that~\eqref{gamma_ass_finite} holds.

Let $u(t)$, $v(t)$ and $w(t)$ be given and consider the functions

\begin{align}\label{V}%
\scalebox{0.92}{\ensuremath{
    \begin{aligned}
        V_u(t) :=& \int_t^T (u(\tau)\T(B\T p(\tau) + r)\!+ x(\tau)\T E\T |B\T p(\tau) + r|)\mathrm d\tau \\
        V_v(t) :=& \int_t^T (v\T(\tau)(H\T p(\tau) -\delta) - x(\tau)\T G\T |H\T p(\tau) -\delta|)\mathrm d\tau \\
        V_w(t) :=& \int_t^T w\T(\tau)(F\T p(\tau) - \gamma)\mathrm d\tau\\
        V(t) :=&\: x(t)\T p(t) + V_u(t) + V_v(t) + V_w(t) 
    \end{aligned}
    }}
\end{align}
where $p$ satisfies \eqref{p_ODE}. Note in \eqref{V} that the integrals vanish when $t=T$, and recall that $p(T) = 0$. It is immediate that $V(T) = x(T)\T p(T) = 0$ and it can be verified that $\dot V(t) = -(s^{\top}x(t) +r^{\top}u(t)-\boldsymbol\gamma^{\top}w(t)-\delta^{\top}v(t))$. Consequently, by the fundamental theorem of calculus, $V(0)$ is the cost in \eqref{minimaxProb_finite} corresponding to the given trajectories $u(t)$, $v(t)$ and $w(t)$.
\\\indent
Let $u\in U(x)$ be given and suppose $(F\T p(t) - \gamma)_i > 0$ for some index $i$ on some nonempty interval $\mc I \subset [0,T]$. Let $v(t) \in \diag(\on{sign}(H\T p(t) + \delta))G x(t)$, which implies $V_v(t)=0$ for all $t$.
Consider the parameter $\varepsilon \in \RR_+$, $\varepsilon >0$ and define 
\begin{align}
    w(t)_j = \begin{cases}
        \varepsilon & \text{if } j = i \text{ and } t\in \mc I\\
        0 & \text{otherwise}
    \end{cases}.
\end{align}
Since $|u| \le Ex$, it follows that
\begin{align*}
    u\T(B\T p + r) + x\T E\T |B\T p + r| \ge 0
\end{align*}
and therefore, by our choice of $v(t)$, we have
\begin{align}\label{V_inequality}
    x_0\T p(0) + V_w(0) \leq V(0).
\end{align}
It follows that $V_w(0)\rightarrow\infty$ as $\varepsilon \rightarrow \infty$, and thus \eqref{V_inequality} implies that the cost in \eqref{minimaxProb_finite} is unbounded for all $x_0$ and all $u\in U(x)$. The optimal cost of \eqref{minimaxProb_finite} is therefore finite only if $F\T p(t) \le \gamma$ for all $t$, and in particular if $F\T p(0) \le \gamma$.

$(i)\Longleftarrow(ii)$ 
Assume that the differential equation~\eqref{p_ODE} has a unique solution and~\eqref{gamma_ass_finite} holds. Let $w=0$. We apply Lemma~\ref{LEMA_cont_minimax_finite} in the Appendix. To apply this Lemma we define
\begin{align}\label{eqs_proofthm1}
    f(t,x,u,v) &:= Ax+Bu+Hv   \\
    g(x,u,\omega)&:= s^{\top}x+r^{\top}u -\delta^{\top}v \notag
\end{align}
and we verify that the Assumptions~\ref{A1},~\ref{A2},~\ref{A3} hold.
It is direct that $f: \mathbb{R}^n \times \mathbb R^m \times \mathbb R^l \longrightarrow \mathbb R^n$, $g: \mathbb{R}^n \times \mathbb R^m \times \mathbb R^l \longrightarrow \mathbb R_+$ are linear and continuous. Moreover, choose 
\begin{align*}
    &\kappa_f= \left\| A\right\|+ \left\|B \right\|\left\|E \right\|+\left\| G\right\|\left\| H\right\|
\end{align*}
and 
\begin{align*}
    \left\|f \right\|&= \left\|Ax+Bu+Hv \right\|\\
    &\leq \left\|A \right\|\left\|x \right\|+\left\| B\right\|\left\| u\right\|+\left\| H\right\| \left\|v \right\|\\
    &\leq \left\|x \right\|\kappa_f \leq (1+\left\|x \right\|)\kappa_f
\end{align*}
Thus, Assumption~\ref{A1} holds. Similarly, it is direct that choosing $\kappa_g= \left\| s\right\|+ \left\|r \right\|\left\|E \right\|+\left\| \delta \right\|\left\| H\right\|$ and $\kappa_{\phi}=0$, Assumption~\ref{A2} holds with $\lambda= (\left\|A \right\|+\|s^{\top}\|)$.
Lastly, Assumption~\ref{A3} follows by the linearity structure of the problem. Therefore, the minimax control problem~\eqref{minimaxProb_finite} is a special case of the general minimax setting~\eqref{gral_minimax_prob_setup_cont_finite} in the Appendix. Next, we verify that the differential equation~\eqref{p_ODE} in $(ii)$ is equivalent to~\eqref{HJB_minimax_cont_finite} in Lemma~\ref{LEMA_cont_minimax_finite}. Recall the cost-to-go function~\eqref{11.2} when $w=0$. Define $J^{*}(t,x)=p^{\top}(t)x$,
\begin{align}\label{aux_HJB_finite}
    &\scalebox{0.95}{$0=\underset{u }{\mathrm{\min} }\hspace{0.2mm}\underset{\omega}{\mathrm{\max}}\left [ g(x,u,\omega)+\nabla_t J^{*}(t,x)+\nabla_x J^{*}(t,x)^{\top}f(x,u,v) \right ] $}\notag\\
    &\scalebox{0.95}{$=\underset{u}{\mathrm{\min} }\hspace{0.3mm}\mathrm{max}_{v}\left [ s^{\top}x +r^{\top}u-\delta^{\top}v+\dot{p}(t)^{\top}x \right. $}\notag \\
    &~~~~~~~~~~~~~~~~~~~~~~~~~~~~~~~~~~~~~~~~~~\scalebox{0.95}{$\left.+p(t)^{\top}(Ax+Bu+Hv) \right ]$} \notag \\
    &\scalebox{0.95}{$=\dot{p}(t)^{\top}x+(s^{\top}+p(t)^{\top}A)x+ \underset{\left | u \right | \leq Ex }{\mathrm{\min} }\left [(r^{\top} + p(t)^{\top}B)u \right ]$} \notag \\
    &~~~~~~~~~~~~~~~~~~~~~~~~~~~~~~~~~~~~~~~\scalebox{0.95}{$+\underset{\left | v \right | \leq Gx}{\mathrm{\max} }\left [(-\delta^{\top}+p(t)^{\top}H )v \right ]$} \notag \\
    &\scalebox{0.95}{$=\dot{p}(t)^{\top}x+(s^{\top}+p(t)^{\top}A)x-\left | r^{\top} + p(t)^{\top}B \right |Ex$}  \\
    &~~~~~~~~~~~~~~~~~~~~~~~~~~~~~~~~~~~~~~~~~~~~~\scalebox{0.95}{$+\!\left | -\delta^{\top} \!+ \! p(t)^{\top}H \right |Gx.$ }\notag
\end{align} 
Given the linear nature of the optimization setting and the policy constraint design, the resulting HJI minimax equation becomes decoupled and the optimizing variables attain their optimal values on the boundary, i.e. $u_i \in \left \{ -E_i x, E_i x \right \}$, $v_i \in \left \{ -G_i x, G_i x \right \}$ for all $i$. Thus, the optimal value of the problem~\eqref{minimaxProb_finite} is given by $J^*(0,x_0)=p^{\top}(0)x_0$ and is achieved by 
\begin{align}\label{u_star}
    \mu(x) &= \mathrm{arg}\min_{ |u| \leq Ex}
    \left [ g(x,u,\omega) + \left (\nabla_x p(t)^{\top}x  \right )(f(x,u,\omega)) \right ]\notag\\
    &= \mathrm{arg}\min_{ |u| \leq Ex} \sum_{i=1}^{m}\left[ \left ( r_{i}^{\top} +p(t)^{\top}B_{i} \right )u_{i} \right]
\end{align}
and

\begin{align}\label{v_star}
    v(t)&=\mathrm{arg}\max_{ |v| \leq Gx}
    \left [ g(x,u,\omega) + \left (\nabla_x p(t)^{\top}x  \right )(f(x,u,\omega)) \right ] \notag\\
    &= \mathrm{arg}\max_{ |v| \leq Gx} \sum_{i=1}^{m}\left[ \left ( -\delta_{i}^{\top} +p(t)^{\top}H_{i} \right )v_{i} \right].
\end{align}

Consider now $w\geq 0$ and let $u=\mu(x)$. It follows from~\eqref{V} that $V_u(0)=0$ and also that the cost for any disturbance is given by 
\begin{align*}
    V(0)=x_0^{\top}p(0)+V_v(0)+V_w(0).
\end{align*}
By Lemma~\ref{lem_monotone},~\eqref{gamma_ass_finite} implies $\boldsymbol{\gamma} \geq F^{\top}p(0) \geq F^{\top}p(t)$ for all $t\ge 0$. Note also that $x$ is increasing in $w$ and that $V_v(0)$ is decreasing in $x$. Thus, the value of $V_v(0)+V_w(0)$ is decreasing in $w$ and its supremum is attained at $w=0$. It follows that~\eqref{u_star} is also optimal for $w \geq 0$, in which case the optimal disturbance is again given by~\eqref{v_star}. This proves that $x_0^{\top}p(0)$, which is finite, is the optimal value of the problem~\eqref{minimaxProb_finite}.

Lastly, because for all $i = 1, \dots, m$ the inequality $\left | u \right |\leq Ex$ restricts $u_{i}$ to the interval $\left [ -E_{i}x, E_{i}x \right ]$, the minimum is attained when the input gradient $r_i+p(t)_i^{\top}B_i$ and $u_{i}$ have opposite signs. 
If the input gradient is zero, then any $u_{i}\in\left [ -E_{i}x, E_{i}x \right ]$ is admissible.
Thus, $u_{i}\in- \mathrm{sign}(r_{i}+p(t)^{\top}B_{i})E_{i}x$. This proves the formula for $K(t)$ in~\eqref{K_cont_finite}.
\end{proof}

For finite time horizon, the optimal feedback policy results in a switching controller, since the control objective changes close to the end of the time interval. Here, the switching behavior is akin to that of a bang-bang controller when the input gradient changes sign, while the sparsity structure of the policy is directly determined by the $E$ matrix in the problem constraints. While switching policies may pose certain limitations in practical implementation, they serve as valuable theoretical tools for evaluating fundamental performance limits and worst-case robustness. Importantly, they can guide the synthesis of more practical approximations by revealing structure and sensitivity in optimal control strategies.

In addition, although the solution to the ODE \eqref{p_ODE} is unique, the same is not necessarily true for the control policy \eqref{K_cont_finite} that achieves the optimal cost. At any time instant for which there occurs a switch in the sign of the input gradient, any choice of the corresponding $u_i\in[-E_i x, E_i x]$ would render the control policy optimal. Naturally, this can be neglected if such behavior occurs on a set of zero measure.
It should be noted, however, that $r_{i}+p(t)^{\top}B_{i}=0$ could also hold on a set with positive measure, meaning that multiple choices of the controller achieve the optimal cost. 
Nevertheless, consistently choosing either $u_i = -E_i x$ or $u_i = E_i x$ when the input gradient is equal to zero yields a bang-bang controller and is sufficient for optimality. 

The next example illustrates how optimal control laws can be non-unique while still yielding trajectories that achieve the same cost value.
\begin{exmp}\label{nonuniqueness_example}
    Let
    \begin{align*}
        A&=\begin{bmatrix}
            -2 & 1 & 0\\
            1 & -2 & 0 \\
            0 & 0 & 1
        \end{bmatrix},& B &=\begin{bmatrix}
            1 & 0\\
            -1 & 0 \\
            0 & 2
        \end{bmatrix},& E&= \begin{bmatrix}
            1 & 1 & 0\\
            0& 0 & 1
        \end{bmatrix}\\
        \\
        F&=H=G=0,& s&=\begin{bmatrix}
            1 & 1 & 1
            \end{bmatrix}^{\top},& r&= \begin{bmatrix}
            0 & 0
            \end{bmatrix}^{\top}.
    \end{align*}
    The solution to the Isaacs equation~\eqref{p_ODE} is
    \begin{align}\label{opt_cost_example}
        p(t)=(1-e^{-t})\begin{bmatrix}
            1&
            1&
            1
        \end{bmatrix}^\top.
    \end{align}
    Then, the input gradient equals
    \begin{align*}
        r^\top + p(t)^\top B = \begin{bmatrix}
            0 & 2 (1-e^{-t})
            \end{bmatrix}
    \end{align*}
    Note that $\text{sign}(2 (1-e^{-t}))=1$ for all $t$. Hence, both $K_1(t)$ and $K_2(t)$
    \begin{align*}
       K_1(t)&=\begin{bmatrix}
           1 & 0\\
           0 & 1
       \end{bmatrix}E,& K_2(t)&=\begin{bmatrix}
           -0.5 & 0\\
           0 & 1
       \end{bmatrix}E, &\forall \hspace{1mm}t\in [0, T]
    \end{align*}
    achieve the same optimal cost~\eqref{opt_cost_example} and are stabilizing, with
    \begin{align*}
        \scalebox{.9}{$A-BK_1(t)=\begin{bmatrix}
            -3&0&0\\
            2&-1&0\\
            0&0&-1
        \end{bmatrix}; \hspace{2mm} A-BK_2(t)=\begin{bmatrix}
            -1&2&0\\
            0&-3&0\\
            0&0&-1
        \end{bmatrix}.$}
    \end{align*}
    Figure~\ref{K1_K2_FIG} illustrates that although both controllers achieve identical performance, their transient responses differ notably. Controller $K_1$ results in more aggressive control action, with a rapid decay in State 1 and an overshoot in State 2. In contrast, controller $K_2$ achieves a more gradual and balanced convergence across all states. This example highlights how different control strategies can meet the same performance criteria while exhibiting different dynamic responses.
    \begin{figure}[h]
        \centering
        \includegraphics[scale=0.75]{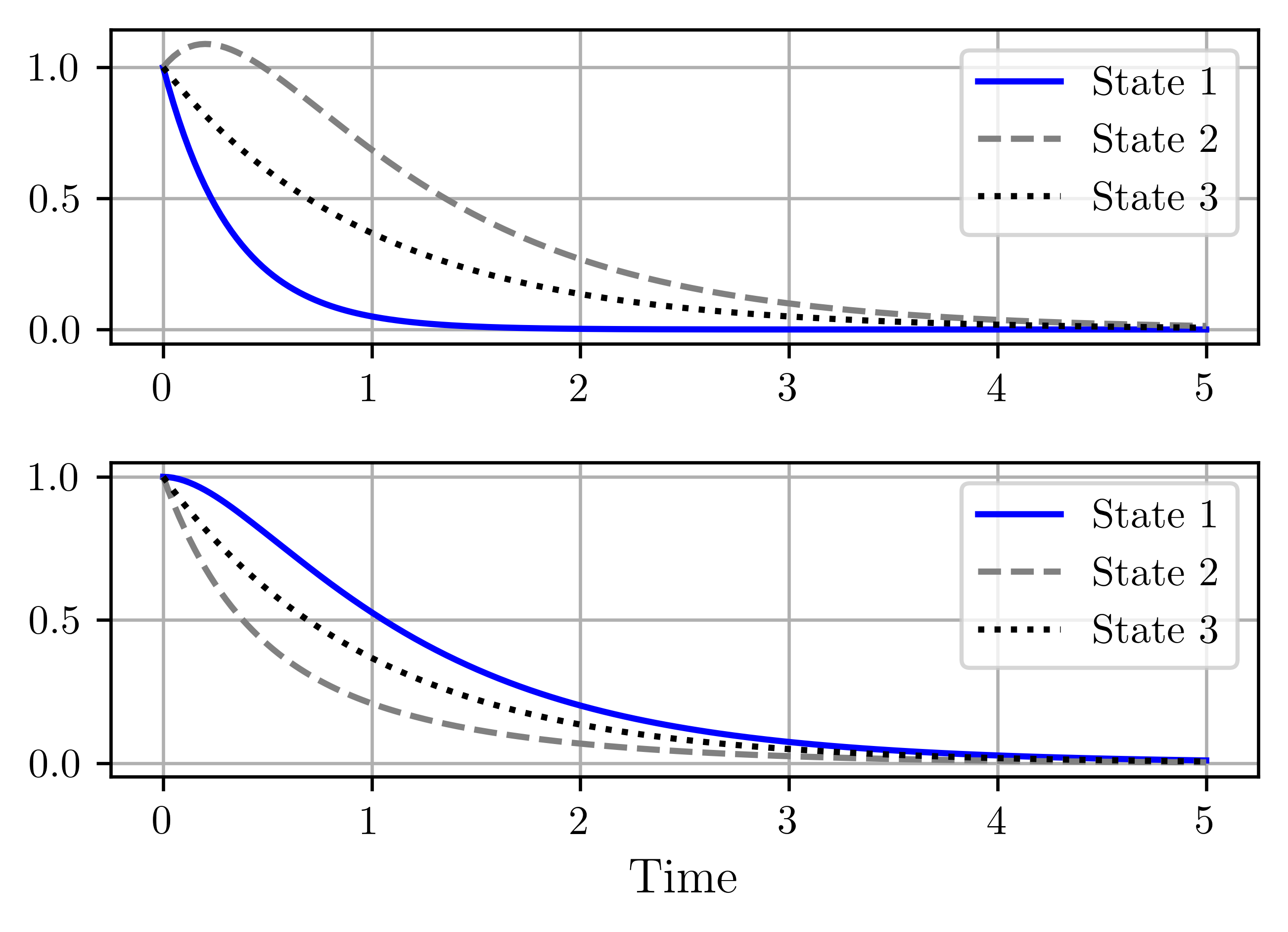}
        \caption{State trajectories of the closed-loop system in Example~\ref{nonuniqueness_example} with initial condition $x_0=\mathbb{1}$, under two different optimal feedback gain matrices (a) $K_1$ and (b) $K_2$.}
        \label{K1_K2_FIG}
    \end{figure}
\end{exmp}

\subsection{Infinite Horizon}\label{subsection_inf_minimax}
Recall the infinite-horizon variant of the optimal control problem specified in \eqref{minimaxProb_finite}.   
\begin{align}\label{minimax_prob_setup_cont_inf}
    &\underset{\mu}{\inf} \hspace{1mm} \underset{w,v}{\sup} \hspace{2mm} \int_{0}^{\infty} \left [ s^{\top}x(\tau) +r^{\top}u(\tau)-\boldsymbol\gamma^{\top}w(\tau)-\delta^{\top}v(\tau)\right ]d\tau \notag \\
            &\mathrm{Subject \hspace{2mm} to} \notag \\
        &~~~~~~~~\dot{x}(t)=Ax(t)+Bu(t)+Fw(t)+Hv(t) \\
        &~~~~~~~~x(0)=x_{0},  \hspace{2mm} x_0 \in \mathbb{R}^n_+ , \hspace{2mm}  u(t)=\mu(x(t)), \notag\\
        &~~~~~~~~\left | u \right | \leq Ex, \hspace{2mm} w \geq 0, \hspace{2mm} \left | v \right | \leq Gx. \notag
\end{align}
Analogous to the finite-horizon case, the following theorem characterizes the optimal solutions to the problem setting~\eqref{minimax_prob_setup_cont_inf}.
\begin{thm}\label{MainThm_minimax_cont_inf}
Let $A \in \mathbb{R}^{n\times n}$, $B = \left [ B_{1}, \dots ,B_{m} \right ] \in \mathbb{R}^{n\times m }$, $F \in \mathbb{R}^{n \times l}_+$, $H \in \mathbb{R}^{n \times c}$, $E = \left [ E_{1}^{\top}, \dots ,E_{m}^{\top} \right ]^{\top} \in \mathbb{R}_{+}^{m \times n}$ such that $E_i^{\top} \neq \mathbb{0}$ for all $i$, $G\in \mathbb{R}^{c \times n}_{+} $, $s \in \mathbb{R}^{n} $, $r \in \mathbb{R}^{m} $, $\boldsymbol\gamma \in \mathbb{R}^{l}_{+}$ and $\delta \in \mathbb{R}^{c}$.  
Suppose that 
\begin{align}\label{ass_a_inf}
    &A -\left | B \right |E -\left | H \right |G  \hspace{1.5mm} \text{Metzler},\\
    &s \geq E ^{\top} \left | r \right | - G^{\top} \left | \delta \right |.
    \label{ass_s_inf}
\end{align}
Then the following statements are equivalent:
\begin{enumerate}[(i)]
\item The optimal control problem~\eqref{minimax_prob_setup_cont_inf} has a finite value for every $x_{0} \in \mathbb{R}_{+}^{n}$;
\item There exists a unique $p\in \mathbb{R}_{+}^{n}$ such that
\begin{align}\label{p_eq_minimax_inf}
    A^{\top}p&= E^{\top}\left | r+B^{\top}p \right |- G^{\top}\left | -\delta +H^{\top}p \right |-s  
\end{align}
and
\begin{align}\label{gamma_ass_inf}
            \boldsymbol{\gamma}  \geq F^{\top}p 
    \end{align}
    which is minimal among all solutions $\hat p$ of \eqref{p_eq_minimax_inf}--\eqref{gamma_ass_inf} in the sense that $p \le \hat p$.
\end{enumerate}
Moreover, if $(i)$--$(ii)$ are satisfied, then 
the infinite horizon \nobreak{minimax} control problem~\eqref{minimax_prob_setup_cont_inf} has minimum value $p^{\top}x_{0}$ and the control law $u(t) = -Kx(t)$ is optimal when 
\begin{align} \label{Kform_minimax_inf_cont}
    K \in \begin{bmatrix}
\mathrm{sign}(r_{1}+ p^{\top}B_{1})E_{1}\\ 
\vdots \\ 
\mathrm{sign}(r_{m}+p^{\top}B_{m})E_{m}
\end{bmatrix}.
\end{align}
\end{thm}

\noindent Note that, Remarks~\ref{rem_K},~\ref{rem_ass_pos},~\ref{rem_1} and~\ref{rem_2} also apply to the infinite horizon case.

Analogous to the finite time horizon case, it may happen that the control law that achieves the optimal cost is not unique. Taking $T\to\infty$ in Example~\ref{nonuniqueness_example} provides an example of such a situation.

\begin{proof}  Let $p_T(t)$ denote the solution to~\eqref{p_ODE} for a fixed end time $T\in \mathbb R_+$. The proof relies on the monotonicity of $p_T(0)$ in $T$ as shown in Lemma~\ref{lem_monotone}, which will be used throughout this proof without further reference.

$(i)\Longrightarrow (ii)$ It is known from Theorem~\ref{theorem_finite} that $J_T(x_0)=x_0^{\top}p_T(0)$ and from Lemma~\ref{lem_monotone} that $p_T(0)$ is monotone nondecreasing in $T$. By assumption
\begin{align}\label{bounded_limit}
    \lim_{T\rightarrow \infty}x_0^{\top}p_T(0) < \infty 
\end{align}
for all $x_0$. Let $p$ be such that $p_i= \lim_{T \rightarrow \infty}e_i^{\top} p_T(0)< \infty$ for all $i=1, \dots, n$, with $e_i$ the $i$-th canonical vector. Then, by linearity, $\lim_{T\rightarrow \infty}x_0^{\top}p_T(0)=x_0^{\top}p$ and the right hand side of the ODE~\eqref{p_ODE} satisfies 
\begin{align*}
    \frac{\partial p_T(0)}{\partial T} \rightarrow 0 \hspace{1mm} \text{as} \hspace{1mm} T\rightarrow \infty
\end{align*}
by monotone convergence. By continuity of the right-hand side of~\eqref{p_ODE}, the algebraic equation~\eqref{p_eq_minimax_inf} holds for $p$.
The inequality~\eqref{bounded_limit} and monotonicity implies that the optimal value of \eqref{minimaxProb_finite} is finite for all $T$. By~\eqref{gamma_ass_finite} this means that $\boldsymbol{\gamma} \ge F^{\top} p_T(0)$. Since $p = \lim_{T\to\infty}p_T(0)$ by construction, this implies the condition~\eqref{gamma_ass_inf}.
By the monotonicity of the dynamics of $p_T(t)$ described by \eqref{p_ODE}, it follows that if any $\hat p \ge 0$ satisfies \eqref{p_eq_minimax_inf} and is an equilibrium of \eqref{p_ODE}, then the trajectory $p_T(t)$ satisfies $p_T(t) \le \hat p$ for all $t$.
Since $p = \lim_{T\to\infty}p_T(0)$ it follows that $p\le \hat p$, and thus $p$ is minimal among all solutions of \eqref{p_eq_minimax_inf}--\eqref{gamma_ass_inf}.

$(i)\Longleftarrow(ii)$
Assume that there exists $p\in \mathbb R^n_+$ which satisfies~\eqref{p_eq_minimax_inf}--\eqref{gamma_ass_inf} and which is minimal among all solutions $\hat p$ of~\eqref{p_eq_minimax_inf},~\eqref{gamma_ass_inf} in the sense that $p \leq \hat p$. 
By the monotonicity of the dynamics of $p_T(t)$ described by \eqref{p_ODE}, it follows that if any $\hat p \ge 0$ satisfies \eqref{p_eq_minimax_inf} and is an equilibrium of \eqref{p_ODE}, then the trajectory $p_T(t)$ satisfies $p_T(t) \le \hat p$ for all $t$.
Since $p_T(0)$ is monotonously nondecreasing, it follows that $\lim_{T\to\infty}p_T(0)\le \hat p$, and thus $p_T(0)$ converges. Let $\bar p := \lim_{T\to\infty}p_T(0)$. It follows that $\bar p$ satisfies \eqref{p_eq_minimax_inf} and that $\bar p \le \hat p$ for all solutions $\bar p$ to \eqref{p_eq_minimax_inf}--\eqref{gamma_ass_inf}. By minimality of $p$ it follows that $p = \bar p$.

Let us denote $\mu_T(x) \in -\diag(\mathrm{sign}(B^{\top}p_T(t)+r))Ex$, which is the optimal feedback of finite time horizon $T$ as described in Theorem~\ref{theorem_finite}.
Consider $\mu(x) \in -\diag(\mathrm{sign}(B^{\top}p+r))Ex$.
Since $p=\lim_{T\to\infty}p_T(0)$ it follows that $\mu_T(x) \to \mu(x)$ as $T\to\infty$.
The corresponding cost for $u=\mu(x)$ is given by
\begin{align}\label{inf_cost_pf}
    &\int_0^{\infty}s^{\top}x+r^{\top}\mu(x)-\delta^{\top}v-\gamma^{\top}w \hspace{1mm} \text d \tau \\
    &=\int_0^{\infty}-\bar p^{\top}(Ax+B\mu(x)+Hv+Fw)\text d \tau+V_v(0)+V_w(0) \notag\\
    &= p^{\top}x_0-\lim_{t \rightarrow \infty} p^{\top}x(t)+ V_v(0)+V_w(0) \notag\\
    &\leq p^{\top}x_0-\lim_{t \rightarrow \infty} p^{\top}x(t)\notag
\end{align}
where $V_u(0)$, $V_v(0)$, $V_w(0)$ are given in~\eqref{V} with $T=\infty$, and where we used that $V_u(0)=0$, $V_v(0) \leq 0$ and $V_w(0) \leq 0$ because $u=\mu(x)$, $|v|\le G x$ and \eqref{gamma_ass_inf}, respectively.
Similarly, the cost for applying $u=\mu_T(x)$ in the finite time horizon problem results in 
\begin{align*}
    p_T(0)^{\top}x_0 + V_{v,T}(0)+V_{w,T}(0).
\end{align*}
Since $\mu_T(x) \to \mu(x)$ as $T\to\infty$, the two costs coincide and thus $\lim_{t \rightarrow \infty} p^{\top}x(t) = 0$. It follows that the optimal cost~\eqref{inf_cost_pf} is $p^{\top}x_0$.%

Finally, the formula for $K$ in~\eqref{Kform_minimax_inf_cont} is derived analogously to the formula for $K(t)$ in the finite horizon case~\eqref{K_cont_finite}. 
\end{proof}

The following theorem proposes an iterative fixed point method to solve the HJI equation~\eqref{p_eq_minimax_inf}. It relies on the iterative method for solving the Isaacs equation of the discrete-time version of the optimal control problem~\eqref{minimax_prob_setup_cont_inf}, derived in~\cite{AlbaEmmaAnders}.
\begin{thm}[Continuous-time value iteration]\label{Cont_VI}
    Let $A\in \mathbb{R}^{n \times n}$, $B\in \mathbb{R}^{n \times m}$, $H \in \mathbb{R}^{n \times c}$, $E \in \mathbb{R}_{+}^{m \times n}$, $s \in \mathbb{R}^{n} $, $r \in \mathbb{R}^{m} $, $\delta \in \mathbb{R}^{c}$ and $h \in \mathbb{R}$. 
    Assume $s \geq E ^{\top} \left | r \right | - G^{\top} \left | \delta \right |$, $A-\left| B\right|E-|H|G$ is Metzler and $A-\left| B\right|E-|H|G + h I_n \geq 0$ with $h>0$. 
    Define $\hat{A}=\frac{1}{h}A + I_n$, $\hat{B}=\frac{1}{h}B$, $\hat{H}=\frac{1}{h}H$, $\hat{E}=E$, $\hat{G}=G$, $\hat{s}=\frac{1}{h}s$, $\hat{r}=\frac{1}{h}r$ and $\hat{\delta}=\frac{1}{h}\delta$. 
    Then, the recursive sequence $\left \{ p_{k} \right \}_{k=0}^{\infty}$ with $p_{0}=0$ and
    \begin{align}\label{p_seq}
        p_{k+1}=\hat{s}+\hat{A}^{\top}p_k-\hat{E}^{\top}\left|\hat{B}^{\top} p_k + \hat{r}\right|+\hat{G}^{\top}\left|\hat{H}^{\top} p_k - \hat{\delta}\right|
    \end{align}
    has a finite limit if and only if there exists $p\in \mathbb{R}^n$ such that~\eqref{p_eq_minimax_inf},
    in which case $p_k \to p$.
\end{thm}
\begin{proof}
    From Theorem 1 in~\cite{AlbaEmmaAnders} it is direct that~\eqref{p_seq} has a finite limit if and only if there exists a $p\in \mathbb{R}^n_+$ such that
    \begin{align}\label{p_disc}
        p=\hat{s}+\hat{A}^{\top}p-\hat{E}^{\top}\left|\hat{B}^{\top} p + \hat{r}\right|+\hat{G}^{\top}\left|\hat{H}^{\top} p - \hat{\delta}\right|.
    \end{align}
    and in this case, $p_k\to p$.
    By definition of $A$, $B$, $H$, $E$, $G$, $s$, $r$ and $\delta$, adding $-p$ to both sides of \eqref{p_disc} and multiplying by $h$ yields
    \begin{align*}
        0&=\hat{s}+(\hat{A}-I_n)^{\top}p-\hat{E}^{\top}\left|\hat{B}^{\top} p + \hat{r}\right|+\hat{G}^{\top}\left|\hat{H}^{\top} p - \hat{\delta}\right| \Rightarrow\\
        0&=h\left(\hat{s}+(\hat{A}-I_n)^{\top}p-\hat{E}^{\top}\left|\hat{B}^{\top} p +\hat{r}\right|+\hat{G}^{\top}\left|\hat{H}^{\top} p - \hat{\delta}\right| \right)\\
        &=s+A^{\top}p-E^{\top}\left|B^{\top} p +r\right|+G^{\top}\left|H^{\top} p -\delta\right|.
    \end{align*}
    Thus, $p$ solves \eqref{p_disc} if and only if $p$ also solves \eqref{p_eq_minimax_inf}.
\end{proof}
\section{The case of positive unconstrained disturbances}\label{section_5}
This section analyzes the explicit solutions of the Isaacs equation for the optimal control problem in the case where the system is affected only by positive, unconstrained disturbances $w$.
Assume therefore for the entire section that $H=G^{\top}=0$. Moreover, within the context of Theorem~\ref{MainThm_minimax_cont_inf} assume also condition~\eqref{gamma_ass_inf}. Then, it follows that the solution $p$ of the minimization problem 
    \begin{align}\label{min_finite_cont}
        &\underset{\mu}{\inf} \hspace{0.3mm}\hspace{1mm}\int_{0}^{\infty} \left [ s^{\top}x(\tau) +r^{\top}u(\tau)\right ]d\tau \notag\\
        &\mathrm{Subject \hspace{2mm} to} \notag \\
        &~~~~~~~~\dot{x}(t)=Ax(t)+Bu(t)\\
        &~~~~~~~~x(0)=x_{0}, \hspace{2mm} u(t)=\mu(x(t)), \notag\\
        &~~~~~~~~\left | u \right | \leq Ex. \notag
    \end{align}
coincides with the solution to the minimax problem~\eqref{minimax_prob_setup_cont_inf}. In this context, the minimal value solution is $p^{\top}x_0$, where $p$ is obtained by solving the algebraic equation
\begin{align}\label{p_min_ODE}
            0&= s + A^{\top}p -E^{\top}\left | r+B^{\top}p \right |. 
        \end{align}
Naturally,  Theorem~\ref{MainThm_minimax_cont_inf} also applies to problem~\eqref{min_finite_cont}. This case is interesting because, as it is shown in this section,~\eqref{p_min_ODE} can be formulated as a linear program.
\subsection{E-Stabilizability and Detectability Conditions}\label{Stab_Detect_subsec}
In the context of~\eqref{min_finite_cont}, two relevant definitions are presented. A pair $(A,B)$ is \nobreak{\textit{E-stabilizable}} if there exists a feedback law $u=-Kx$ with  $\left|u \right|\leq Ex$ such that $A-BK$ is Hurwitz. The system is \textit{detectable} if all unobservable states are stable. A detectability characterization for positive systems can be stated as follows: %
\begin{prop}[Proposition 3~\cite{math_detect}]\label{Prop_math}
Let $A \in \mathbb R^{n \times n}$, $B \in \mathbb R^{n \times m}$ and $C \in \mathbb R^n$. Consider the autonomous system
\begin{align*}
        \dot{x}=Ax,\qquad
        y=Cx
\end{align*}
\noindent where $A$ is a Metzler matrix and $C\geq 0$. The pair $(C,A)$ is detectable if and only if $Cv>0$ for any nonnegative eigenvector $v$ corresponding to a nonnegative eigenvalue $\lambda$ of \nobreak{$A$}.
\end{prop}
The next theorem presents an \textit{a posteriori} condition that determines if an optimal control law to \eqref{minimaxProb_finite} stabilizes the system dynamics.
\begin{thm}\label{bellman_detectable}
    Suppose that the Isaacs equation \eqref{p_min_ODE} has a finite \nobreak{solution} $p$. Let $K$ satisfy \eqref{Kform_minimax_inf_cont}.
    If the pair $(s\T-r\T K, A-BK)$ is detectable, then $u(t) = - Kx(t)$ stabilizes the system.
\end{thm}
\begin{proof}
    Since $p$ solves the Isaacs equation, the controller $K$ achieves the optimal cost.
    The optimal cost in closed-loop equals $J^*(x_0) = \int_0^t(s\T-r\T K)x(\tau)\d \tau$ and is optimal. It satisfies $J^*(x_0) = p \T x_0$ and is bounded. Therefore, $(s\T-r\T K)x(\tau) \to 0$ as $\tau\to\infty$. If the pair $(s\T-r\T K, A-BK)$ is detectable, this implies that the closed-loop matrix $A-BK$ is Hurwitz.
\end{proof}

Detectability of the pair $(s\T-r\T K, A-BK)$ can be verified through Proposition~\ref{Prop_math}, since the closed-loop system is again a positive system and $s\T-r\T K$ is nonnegative thanks to~\eqref{ass_s_inf}. The next corollary gives a simple but restrictive \textit{a priori} condition that implies observability for all of the controllers that can be produced by Theorem~\ref{MainThm_minimax_cont_inf}. A more sophisticated condition is introduced in Theorem~\ref{thm_detect_cond}.
\begin{cor}\label{cor_suff_cond_detec}
    Suppose $s - E^{\top} \left|r \right| > 0$, then $(s\T-r\T K, A-BK)$ is detectable for any $K$ satisfying $|K|\le E$.
    Moreover, if $(i)$--$(ii)$ in Theorem~\ref{MainThm_minimax_cont_inf} hold with $G=H^{\top}=0$ and $K$ that satisfies \eqref{K_cont_finite}, then $u(t) = - K x(t)$ stabilizes the system.
\end{cor}
\begin{proof}
    Since $s - K^{\top} r \ge s - E^{\top} \left|r \right| > 0$ we have $(s - K^{\top} r)\T v > 0$ for every nonzero nonnegative eigenvector $v$ of $A-BK$. Hence $(s\T-r\T K, A-BK)$ is detectable by Proposition~\ref{Prop_math}. The remainder follows from Theorem~\ref{bellman_detectable}.
\end{proof}

In the next theorem, we demonstrate that the detectability of $(s^{\top}-K^{\top}r, A-BK)$ is guaranteed if every nonzero state incurs a positive cost regardless of the policy applied.
\begin{thm}\label{thm_p_pos_detect}
    Assume conditions $(i)$ and $(ii)$ in Theorem~\ref{MainThm_minimax_cont_inf} hold with $G=H^{\top}=0$. Let $K$ satisfy~\eqref{Kform_minimax_inf_cont}. If there exists a solution $p \in \mathbb{R}^n_+$ to~\eqref{p_min_ODE} such that $p > 0$ and the spectral abscissa $\alpha(A-BK) \neq 0$, then the pair $(s^{\top}-r^{\top}K, A-BK)$ is detectable and $u(t)=-Kx(t)$ is a stabilizing policy. Moreover, if $\alpha(A-BK)=0$ then the pair ($s^{\top}-r^{\top}K$, $A-BK$) is not detectable.
\end{thm}

\begin{proof}
    Suppose there exists a vector $p \in \mathbb{R}^n_+$ solving~\eqref{p_min_ODE} with $p > 0$. Note that 
    \begin{align*}
        K^{\top}(B^{\top}p + r) = E^{\top}|B^{\top}p + r|.
    \end{align*}
    Substituting this into~\eqref{p_min_ODE}, it follows that
    \begin{align}\label{aux_detect_cond}
        -(A-BK)^{\top}p = s - K^{\top}r.
    \end{align} 
    
    If $A-BK$ is Hurwitz, then all unobservable states of $A-BK$ are asymptotically stable. Thus, the pair $(s^{\top}-r^{\top}K, A-BK)$ is detectable and $u(t)=-Kx(t)$ is a stabilizing policy.

    Suppose $A-BK$ is not Hurwitz. Since $A-|B|E$ is Metzler, $A-BK$ is also Metzler. Therefore, $A-BK$ has a nonnegative Perron eigenvector $v_p \neq 0$ associated with the Perron root $\lambda_p=\alpha(A-BK)\neq 0$ satisfying $(A-BK)v_p=\lambda_p v_p.$ 
    The right-hand side of~\eqref{aux_detect_cond} can be rewritten as
    \begin{align}\label{aux_detect_perron}
        p^{\top}(-(A-BK))v_{\mathrm{p}} = p^{\top} (-\lambda_{\mathrm{p}}) v_{\mathrm{p}}. 
    \end{align}
    Since $p>0$, $v_{\mathrm{p}} \geq 0$ and $\lambda_{\mathrm{p}} > 0$ it follows that $p^{\top} (-\lambda_{\mathrm{p}}) v_{\mathrm{p}}  <0$. 
    This contradicts the left-hand side of~\eqref{aux_detect_cond}, $(s - K^{\top}r)^{\top}v_{\mathrm{p}} \geq 0$. Therefore, if $p>0$ and $\lambda_{\mathrm{p}} \neq 0$ the inequality cannot hold. Consequently, if there exists a solution $p>0$ to~\eqref{p_min_ODE} and $\alpha(A-BK)\neq 0$ the pair $(s^{\top}-r^{\top}K, A-BK)$ is detectable and $u(t)=-Kx(t)$ is a stabilizing policy.
    
    Finally, if $\alpha(A-BK)=0$, the Perron root is $\lambda_{\mathrm{p}}=0$ with $v_{\mathrm{p}} \geq 0$. Therefore,~\eqref{aux_detect_perron} gives $p^{\top} (-\lambda_{\mathrm{p}}) v_{\mathrm{p}}=(s-K^{\top}r)^{\top}v_{\mathrm{p}} =0$. Thus, $(s\T-r\T K, A-BK)$ is not detectable by Proposition~\ref{Prop_math}.
\end{proof}

Recall the observability condition in~\cite[Assumption 4.1]{yuchao}.
\begin{align}\label{obs_cond_disc}
    (s^{\top}-|r^{\top}|E)\sum_{i=0}^{n-1}(A-|B|E)^i>0.
\end{align}
Condition~\eqref{obs_cond_disc} is based on the sum of the observability matrix in the discrete-time framework. In continuous time, the observability matrix involves the integral $\int_{0}^{\infty}e^{(A-|B|E)t}dt$. This integral does not converge for the eigenvalues $\lambda$ of $A-|B|E$ for which $\mathrm{Re}(\lambda)\geq 0$. Thus, the observability condition
\begin{align*}
     \int_{0}^{\infty}e^{(A-|B|E)^{\top}}(s-E^{\top}|r|)dt>0
\end{align*}
ensures that no eigenvector $v \neq 0$ associated with an unstable eigenvalue $\lambda$ such that $\mathrm{Re}(\lambda)\geq 0$, %
is unobservable.
Inspired by~\eqref{obs_cond_disc}, the next theorem proposes an \textit{a priori} detectability condition for the continuous-time setting.
\begin{thm}\label{thm_detect_cond}
    Assume $(i)$--$(ii)$ in Theorem~\ref{theorem_finite} and~\ref{MainThm_minimax_cont_inf} hold with $G=H^{\top}=0$. Let $h \in \mathbb{R}$ such that  $A-\left| B\right|E + h I_n \geq 0$ and $h>0$. Suppose also that 
    \begin{align}\label{obs_cond_cont}
         (s^{\top}-|r^{\top}|E) \sum_{i=0}^{n-1}\left(\frac{1}{h}(A-|B|E)+I \right)^i>0.
     \end{align}
    Then, $(s\T-r\T K, A-BK)$ is detectable for any $K$ satisfying $|K|\le E$.
\end{thm} 
\begin{proof}
Suppose~\eqref{obs_cond_cont} holds. 
Let $M=A-BK$, $N=s-K^{\top}r$, $\hat M=A-|B|E$ and $\hat N=s-E^{\top}|r|$. Note that $M\ge \hat M$ and $N \ge \hat N$. 
The ODE~\eqref{p_ODE} can be rewritten as $-\dot{p}(t)=Mp(t)+N$ with $p(T)=0$. Define $w(t)=p(T-t)$ and suppose $T>0$, then $w(0)=0$ and $w(t)$ for $t\in[0,T]$ satisfies 
\begin{align}\label{ODE_sol_w}
    w(t)=\int_0^{t}e^{M^{\top} \tau}N d\tau.
\end{align}
It is possible to approximate the integral by
\begin{align}\label{approx_exp}
 w(t)=\int_0^{t}e^{M^{\top}\tau}Nd \tau &= \lim_{h\rightarrow \infty} \frac{1}{h} \sum_{i=0}^{\left\lfloor \frac{t}{h} \right\rfloor -1}( \frac{1}{h}M^{\top} +I)^i N
\end{align}
Because $M + hI_n \geq \hat M + hI_n \geq 0$ and $N \ge \hat N\ge 0$, the summation in~\eqref{approx_exp} can be lower bounded by
\begin{align*}
    \lim_{h\rightarrow \infty} \frac{1}{h} \sum_{i=0}^{n-1}\left( \frac{1}{h}\hat M^{\top} + I \right)^i \hat N > 0
\end{align*}
which is positive by assumption.
For a nonnegative vector ${v\neq 0}$ corresponding to an eigenvalue $\lambda \geq 0$ of $M$ it holds that
\begin{align*}
    v^{\top}w(t) &= v^{\top}\lim_{h\rightarrow \infty} \frac{1}{h} \sum_{i=0}^{\left\lfloor \frac{t}{h} \right\rfloor -1}(\frac{1}{h}M^{\top} +I)^i N\\
    & = \lim_{h\rightarrow \infty} \frac{1}{h} \sum_{i=0}^{\left\lfloor \frac{t}{h} \right\rfloor -1}(\frac{\lambda}{h} + 1)^i v^{\top} N = \big(\int_0^t e^{\lambda\tau}d\tau\big) v^{\top} N.
\end{align*}
Clearly $\int_0^t e^{\lambda\tau}d\tau > 0$. At the same time
\begin{align*}%
    v^{\top}w(t) &= v^{\top}\lim_{h\rightarrow \infty} \frac{1}{h} \sum_{i=0}^{\left\lfloor \frac{t}{h} \right\rfloor -1}( \frac{1}{h}M^{\top} +I)^i N\\
    &\ge v^{\top}\lim_{h\rightarrow \infty} \frac{1}{h} \sum_{i=0}^{n -1}( \frac{1}{h}\hat M^{\top} +I)^i \hat N > 0.
\end{align*}
Thus, $v^{\top} N = v^{\top} (s-K\T r ) >0$ and $(s\T-r\T K, A-BK)$ is detectable by Proposition~\ref{Prop_math}.
\end{proof} 
\subsection{Linear Programming}\label{LP_subsect}
It is possible to obtain a solution to the Isaacs equation \eqref{p_min_ODE} through the linear program 
\begin{align}\label{LP}
        &\mathrm{Maximize} \hspace{2mm} \mathbb{1}^{\top}p \hspace{1mm} \mathrm{over} \hspace{1mm} p \in \mathbb{R}^{n}_{+}, \hspace{1mm} \zeta \in \mathbb{R}^{m}_{+} \notag\\
        &\mathrm{Subject} \hspace{1mm} \mathrm{to} \hspace{2mm} A^{\top}p \geq E^{\top}\zeta -s   \\
        &\hspace{15mm} -\zeta \leq r + B^{\top}p \leq \zeta \notag . %
    \end{align}
Lemma~\ref{lem_thm_LP} characterizes the boundedness of  \eqref{LP} by considering the dual linear program.
\begin{lem}\label{lem_thm_LP}
    The following are equivalent:
    \begin{enumerate}[(i)]
        \item The primal linear program~\eqref{LP} has a bounded solution;
        \item The dual linear program
        \begin{align}\label{dual_LP}
        &\mathrm{Minimize}\quad  
            s\T x + r\T u~\text{ with }x \in \mathbb{R}^{n}_{+}
            \notag \\
            &\mathrm{Subject \hspace{1mm} to}\quad 
             Ax + Bu \le - \mb 1
            \\
            &\hspace{14mm}  -Ex \le u \le Ex\notag
        \end{align} %
        is feasible;
        \item There exists a $D\in\RR^{m\times m}$ satisfying $-I\le D \le I$ such that $A-BDE$ is Hurwitz.
    \end{enumerate}
\end{lem}
\begin{proof}
    ($i) \Longleftrightarrow (ii$): This follows from the weak duality of the linear programs, and the fact the the primal problem is feasible. Indeed, taking $p=0$ and $\zeta = |r|$ renders the primal linear program feasible since $s \ge E\T |r|$.

    ($ii) \Longrightarrow (iii$): 
    For any feasible $(x,u)$ let $D$ such that $-I\le D \le I$ and $u = -DEx$.
    It follows that
        $(A-BDE)x \le -\mb 1 < \mb 0.$
    The matrix $A-BDE$ is Metzler since $A-|B|E$ is Metzler and $-I\le D \le I$. 
    Since $x\ge 0$, it follows by item 1 and 18 of \cite[Thm. 5.1]{fiedler1986} that $A-BDE$ is Hurwitz. 

    ($iii) \Longrightarrow (ii$): 
    By item 2 and 18 of \cite[Thm. 5.1]{fiedler1986} there exists a $v>0$ such that $(A-BDE)v < \mb 0$. Let $\alpha = \min_i |(A-BDE)_i|$ and take $x = \frac 1 \alpha v$ and $u=-DEx$, then $Ax + Bu = (A-BDE)x \le -\mb 1$, $-Ex \le u \le Ex$ and $x\ge \mb 0$.
\end{proof}
\begin{thm}\label{thm_LP}
    If the linear program \eqref{LP} has a bounded value, then the optimizer $p$ solves the Isaacs equation \eqref{p_min_ODE}.
\end{thm}
\begin{proof}
    Let $p$ and $\zeta$ optimize \eqref{LP}, then $A\T p = E\T \zeta - s$ and $\zeta = |r+B\T p|$, and thus \eqref{p_eq_minimax_inf} with $G=H\T=0$ is satisfied.
\end{proof}

Note that the existence of a bounded solution to the LP~\eqref{LP} is not necessary for the existence of solution to the HJI equation. We illustrate this case by an example:
\begin{exmp}
    Let
    \begin{align*}
        A&=\begin{bmatrix}
            -2 & 1 & 0\\
            1 & -2 & 0 \\
            0 & 0 & 1
        \end{bmatrix}; \hspace{1mm} B =\begin{bmatrix}
            1 \\
            -1 \\
            0
        \end{bmatrix}; \hspace{2mm} E^{\top}= \begin{bmatrix}
            1 \\ 1 \\ 0
        \end{bmatrix}\\
        s&=\begin{bmatrix}
            1 & 1 & 0
            \end{bmatrix}^{\top}; \hspace{2mm} r= 0.
    \end{align*}
    The third state is unstable and not detectable (nor stabilizable).
    The solution to the Isaacs equation~\eqref{p_min_ODE} is 
    \begin{align}
        p=\begin{bmatrix}
            1 & 1 & 0
        \end{bmatrix}^{\top}.
    \end{align}
    However, the third entry of the $p$ vector in the linear program is unbounded. 
\end{exmp}
Moreover, a system that is not detectable might still be stabilizable, but the Isaacs equation might fail to identify the stabilizing solution. This is illustrated in the next example: 
\begin{exmp}
    Let
    \begin{align*}
        A&= \begin{bmatrix}
            \frac{1}{2} & 1 & 0\\
            0 & -1 & 0 \\
            0 & 0 & -1
        \end{bmatrix}; \hspace{1mm} B= \begin{bmatrix}
            1 \\
            0 \\
            0
        \end{bmatrix}; \hspace{1mm} E^{\top} = \begin{bmatrix}
            1 \\
            0 \\
            0
        \end{bmatrix}\\
        s&= \begin{bmatrix}
            0 & 2 & 3
        \end{bmatrix}^{\top}; \hspace{1mm} r=0.
    \end{align*}
    The Perron Frobenius eigenvalue of $A$ is $\lambda_p = 0.5$ indicating that the first state is unstable. Furthermore, the detectability condition fails since $$(s^{\top}-|r^{\top}|E)v_p= \begin{bmatrix}
        0 & 2 & 3
    \end{bmatrix} \begin{bmatrix}
        1 & 0 & 0
    \end{bmatrix}^{\top}=0$$ where $v_p$ is the eigenvector associated with $\lambda_p$. This implies that the first state is not detectable. However, the system is stabilizable: Choosing $K=E$ yields
    \begin{align*}
       A-BK= \begin{bmatrix}
            -\frac{1}{2} & 1 & 0\\
            0 & -1 & 0 \\
            0 & 0 & -1
        \end{bmatrix}, \hspace{2mm} \lambda_{\mathrm{p}}= -0.5.
    \end{align*}
In this context, the solution to the Isaacs equation~\eqref{p_min_ODE} is
\begin{align}
    p= \begin{bmatrix}
        0 & 2 & 3
    \end{bmatrix}^{\top}
\end{align}
with $K=\mathrm{sign}(0)E$, i.e. any $K \in [-E, E]$ is minimizing. However, $K=-E$ does not give a stable closed-loop, so not all minimizing solutions are stabilizing.
\end{exmp}
Lemma~\ref{lem_thm_LP} implies that the linear program~\eqref{LP} converges only if $(A, B)$ stabilizable by a feedback $u=-Kx$ with $\left|u \right|\leq Ex.$ Therefore, a solution to the Isaacs equation~\eqref{p_min_ODE} is an optimizer to~\eqref{LP} as long as the closed-loop system is detectable.

Next, a partial converse to Theorem~\ref{thm_LP} is presented:
\begin{thm}
    If $p$ solves the Isaacs equation \eqref{p_eq_minimax_inf}, $K$ satisfies \eqref{Kform_minimax_inf_cont}, and $K$ is such that
    \begin{align}\label{obs_cond}
        \left(s-K^{\top} r, A-BK \right)
    \end{align}
    is detectable, then $p$ maximizes \eqref{LP}, with $\zeta = |r + B^\top p|$.
\end{thm}
\begin{proof}
    For the primal linear program to be bounded we require by Lemma~\ref{lem_thm_LP} that there exists a feedback law $u(t)=-Kx(t)$ with $K=DE$ and $|D|\le I$ that stabilizes the system.
    By Theorem~\ref{bellman_detectable}, it follows that any $K$ that satisfies \eqref{Kform_minimax_inf_cont} stabilizes the system.
    Taking $p$ in \eqref{LP} as the solution to \eqref{p_eq_minimax_inf} and $\zeta = |r + B\T p|$ satisfies constraints of \eqref{LP}.
    If $p$ would not be optimal, then applying the iteration \eqref{Cont_VI} to $p$  would result in a $\tilde p$ for which $\mathbb 1\T p \le \mathbb 1\T \tilde p$ \cite{yuchao}. However, since $p$ solves the Isaacs equation we have $p=\tilde p$, and thus $p$ is optimal.
\end{proof}

From the proof of Lemma~\ref{lem_thm_LP} it is clear that the dual linear program \eqref{dual_LP} always generates a stabilizing controller. We use this fact to show that the primal linear program \eqref{LP} also generates at least one stabilizing controller.
\begin{thm}
    If the primal linear program~\eqref{LP} has a bounded solution $p$, then all controllers $u(t)=-Kx(t)$ satisfying
    \begin{align*}
        K \in \diag(\operatorname{sign}(B\T p + r))E
    \end{align*}
    achieve the optimal cost and at least one stabilizes the system, \text{i.e.}, the resulting closed-loop matrix $A-BK$ is Hurwitz.
    Moreover, if $(s\T-r\T K, A-BK)$ is detectable, then $K$ always stabilizes the system.
\end{thm}
\begin{proof}
    The primal and dual linear programs~\eqref{LP} and \eqref{dual_LP} achieve the same cost due to strong duality (\textit{e.g.}, see \cite[Ch.~2]{dantzig2003}). Let $p^*$, $\zeta^*$, $x^*$ and $u^*$ be such as the optimal cost of both programs is achieved. That is, $s\T x^* + r\T u^*=\mb 1\T p^*$. Thanks to optimality of the constraints we have
    \begin{align*}
        &|r + B\T p^*| = \zeta^*,
        \quad
        E\T\zeta^* = s + A\T p^*,\\
        &Ax^* + Bu^* = -\mb 1,
        \quad
        |u| = Ex^*.
    \end{align*}
    Let $D=\diag(\operatorname{sign}(u^*))$ such that $u^* = DEx^*$.
    If $x_i^* = 0$, note that $((A+BDE)x^*)_i \ge 0$ since all off-diagonal entries of $A+BDE$ are nonnegative. This violates $(A+BDE)x^* = - \mb 1$ and therefore  $x^*>0$.
    By assumption, $E$ is nonnegative and has no all-zero rows, and so $Ex^* > \mb 0$.
    The established equalities yield
    \begin{align*}
        |r + B\T p^*|\T E x^* 
        = \zeta^*\mathstrut\T E x^* 
        = s\T x^* + p^*\mathstrut\T Ax^*
        \\
        = s\T x^* - p^*\mathstrut\T \mb 1 - p^*\mathstrut\T Bu^*
        = -r\T u^* - p^*\mathstrut\T Bu^*
        \\
        = -(r + B\T p^*)\T u^*
        = -(r + B\T p^*)\T D E x^*.
    \end{align*}
    It follows that equality holds if $D \in -\diag(\operatorname{sign}(r + B\T p^*))$.
    The controller is not necessarily unique if there exists an index $i$ such that $(r + B\T p^*)_i = 0$.
    At least one of such $D$ satisfies $D = \diag(\operatorname{sign}(u^*))$ and therefore stabilizes the system.

    The optimal cost $J^*(x_0) = \int_0^t(s\T-r\T K)x(\tau)\d \tau$ satisfies $J^*(x_0) = p^*\mathstrut \T x_0$ and is bounded. Therefore, $(s\T-r\T K)x(\tau) \to 0$ as $\tau\to\infty$. If the pair $(s\T-r\T K, A-BK)$ is detectable, this implies that $A-BK$ is Hurwitz.
\end{proof}

\section{\texorpdfstring{$L_{1}-$} -induced gain Minimization Problem}\label{section_l1gain}

Recall that the $L_1$-induced gain of a system is the maximum ratio of the $L_1$ norm of the system's output to the $L_1$ norm of the control and disturbance input. Thus, it measures the maximum amplification of the input disturbances and control signals to the system's output. Formally, we define the $L_{1}-$induced gain of the system~\eqref{dynamics} with respect to a disturbance $\omega$ as 
\begin{align}\label{l1_def}
        \left \| G_{\mu, \omega} \right \|_{1-\mathrm{ind}}=\sup_{\omega \neq 0}\frac{\left\| z \right\|_{1}}{\left\| \omega \right\|_{1}}=\sup_{\omega \neq 0}\frac{\int_{0}^{\infty}\left[s^{\top}x(\tau) + r^{\top}u(\tau)\right] \hspace{1mm} d\tau}{\int_{0}^{\infty} \hspace{1mm} \mathbb{1}^{\top}\omega(\tau) \hspace{1mm} d\tau}.
\end{align}
In this section, we derive the $L_1$-induced gain of the closed-loop system dynamics~\eqref{dynamics} with respect to the disturbance $w\geq 0$ and the worst case disturbance of~\eqref{min_finite_cont} and~\eqref{minimax_prob_setup_cont_inf}. The same results are obtained for the infinite time setting in ~\cite[Sec. 6]{BLANCHINI_REV2}, but the current paper extends this analysis by characterizing the $L_1$-induced gain of the system.
In addition we give tight bounds on the disturbance penalty to ensure that the optimal control problem achieves a finite cost.
In this section we motivate the bounds \eqref{gamma_ass_finite} and \eqref{gamma_ass_inf} in the context of the $L_1$-induced gain of the disturbance.
Assume in Theorems~\ref{theorem_finite} and~\ref{MainThm_minimax_cont_inf} that $G = H\T = 0$ and $w \neq 0$. Let $\boldsymbol\gamma = \gamma \hspace{0.5mm} \mathbb{1}$. Let the $L_1-$induced gain of the system in Figure~\ref{diagram} with respect to the disturbance $w$ be bounded by a parameter $\gamma^*$, 
    \begin{align}\label{l1_w}
        \left \| G_{\mu, w} \right \|_{1-\mathrm{ind}}=\sup_{w \geq 0}\frac{\left\| z \right\|_{1}}{\left\| w \right\|_{1}} \leq \gamma^*.
    \end{align}
It is immediate from~\eqref{l1_def} that~\eqref{l1_w} gives
    \begin{align*} &\sup_{w \geq 0}\int_{0}^{T}\left[s^{\top}x(\tau) + r^{\top}u(\tau) -\gamma^* \mathbb{1}^{\top} w(\tau)\right] \hspace{1mm} d\tau \leq 0.
    \end{align*} 
    Therefore, given a fixed time horizon $T \in \left[0, \infty \right]$
    \begin{align}\label{cost_l1_gain}
    \mathrm{arg} \hspace{1mm}& \max_{w \geq 0} \left[\int_{0}^{T}\left[s^{\top}x(\tau) + r^{\top}u(\tau) -\gamma \mathbb{1}^{\top}w(\tau)\right] \hspace{1mm} d\tau \right] \notag \\
    &=  \left\{\begin{matrix}
    0 \hspace{3mm} \gamma \geq \gamma^{*}\\
    \infty \hspace{3mm} \gamma < \gamma^{*}
    \end{matrix}\right.  
    \end{align}
    holds. Thus, the cost function~\eqref{cost_l1_gain} has a finite value when $\gamma \geq \gamma^*$.
    From Theorems~\ref{theorem_finite} and~\ref{MainThm_minimax_cont_inf} we know that there exists a finite solution to the optimal control problem~\eqref{minimaxProb_finite} and~\eqref{minimax_prob_setup_cont_inf} if and only if the conditions in \eqref{gamma_ass_finite} and \eqref{gamma_ass_inf} hold. Thus, the minimum $L_1$-induced gain of the system in Figure~\ref{diagram} with $G=H\T=0$ is
    \begin{align}
        \boldsymbol{\gamma^{*}} = \boldsymbol{\gamma}  =
        \scalebox{.9}{\ensuremath{\begin{cases}
              F^{\top}p(0).&\text{if } T < \infty\\
            F^{\top}p &\text{if } T = \infty.
        \end{cases}}}
    \end{align}
The worst-case disturbance can then be found in~\cite{BLANCHINI_REV2}.

\section{Example: Water-flow Network}\label{sect_examples}
This section introduces a model of a water-flow network to which the presented framework is applied. It is illustrated that the technical a priori requirement~\eqref{ass_a} of Theorem~\ref{theorem_finite} and~\eqref{ass_a_inf} of Theorem~\ref{MainThm_minimax_cont_inf}, which assumes that $A-|B|E-|H|G$ is Metzler, can be relaxed to $A-|B|E$ when the worst-case disturbance policy always gives a closed-loop matrix that is Metzler. We demonstrate that the optimal control policy can counteract the effect of the optimal disturbance policy and can compensate an a priori overestimation of the disturbances.

Consider a river system with a downstream flow represented by $\beta$ and the dissipation capacity $\alpha$. The river flow is segmented into $n$ sections, where the state $x$ reflects the volume of water in each section, with initial volume one. Dams positioned at these sections act as controllers that modulate the downstream flow, i.e. the control signal $u$ models the water volume that is restricted from flowing downstream. In this context, disturbances $v$ arise due to leakage effects, where water is dissipated more than is accounted for, resulting in 
extraneous flow towards downstream reservoirs. These leakages intensify as the disturbances propagate downriver, influencing the water volumes of subsequent sections. Additionally, we denote by
$w$ an unconstrained positive disturbance in the form of rain.

\begin{figure}[h]
\centering
\tikzset{every picture/.style={line width=0.75pt}} %

\begin{tikzpicture}[x=0.72pt,y=0.72pt,yscale=-0.9,xscale=0.9]
\draw    (100,101) -- (153,101) -- (159.67,101) ;
\draw    (159.67,80.33) -- (159.67,101) ;
\draw    (100,80.33) -- (100,101) ;
\draw    (189.33,130.33) -- (242.33,130.33) -- (249,130.33) ;
\draw    (249,109.67) -- (249,130.33) ;
\draw    (189.33,109.67) -- (189.33,130.33) ;
\draw    (370.67,191) -- (423.67,191) -- (430.33,191) ;
\draw    (430.33,170.33) -- (430.33,191) ;
\draw    (370.67,170.33) -- (370.67,191) ;
\draw   (181.67,99.67) -- (178.33,99.67) -- (178.33,122.33) -- (181.67,122.33) -- cycle ;
\draw   (362.33,161.67) -- (359,161.67) -- (359,184.33) -- (362.33,184.33) -- cycle ;
\draw    (139,84.33) .. controls (131.08,78.39) and (186.54,46.97) .. (207.7,96.8) ;
\draw [shift={(208.33,98.33)}, rotate = 248.33] [color={rgb, 255:red, 0; green, 0; blue, 0 }  ][line width=0.75]    (10.93,-3.29) .. controls (6.95,-1.4) and (3.31,-0.3) .. (0,0) .. controls (3.31,0.3) and (6.95,1.4) .. (10.93,3.29)   ;
\draw  [line width=0.75] [line join = round][line cap = round] (100.33,92) .. controls (102.66,92) and (114,95.39) .. (115.67,94) .. controls (116.87,92.99) and (114.79,89.6) .. (116.33,89.33) .. controls (121.39,88.45) and (126.58,88.66) .. (131.67,88) .. controls (131.8,87.98) and (136.04,84.02) .. (136.33,84) .. controls (148.19,83.15) and (152.79,89.33) .. (160.33,89.33) ;
\draw  [line width=0.75] [line join = round][line cap = round] (190.33,122.67) .. controls (193.61,122.67) and (193.57,120.14) .. (194.33,120) .. controls (198.9,119.17) and (203.58,122.92) .. (206.33,122) .. controls (208.86,121.16) and (210.34,117.33) .. (211.67,116.67) .. controls (214.14,115.43) and (219.09,115.15) .. (222.33,115.33) .. controls (230.4,115.78) and (240.77,124.12) .. (249,120) ;
\draw    (332.33,155) .. controls (324.41,149.06) and (379.87,117.64) .. (401.04,167.47) ;
\draw [shift={(401.67,169)}, rotate = 248.33] [color={rgb, 255:red, 0; green, 0; blue, 0 }  ][line width=0.75]    (10.93,-3.29) .. controls (6.95,-1.4) and (3.31,-0.3) .. (0,0) .. controls (3.31,0.3) and (6.95,1.4) .. (10.93,3.29)   ;
\draw  [line width=0.75] [line join = round][line cap = round] (371.67,181.33) .. controls (374.06,180.14) and (379.98,176.64) .. (383,177.33) .. controls (384.39,177.65) and (384.32,180.22) .. (385.67,180.67) .. controls (392.66,183) and (399.39,180) .. (406.33,180) .. controls (409.45,180) and (412.58,181.13) .. (415.67,180.67) .. controls (420.41,179.96) and (428.15,181.88) .. (429.67,177.33) ;
\draw    (221.67,110.33) .. controls (213.75,104.39) and (269.21,72.97) .. (290.37,122.8) ;
\draw [shift={(291,124.33)}, rotate = 248.33] [color={rgb, 255:red, 0; green, 0; blue, 0 }  ][line width=0.75]    (10.93,-3.29) .. controls (6.95,-1.4) and (3.31,-0.3) .. (0,0) .. controls (3.31,0.3) and (6.95,1.4) .. (10.93,3.29)   ;
\draw  [dash pattern={on 0.84pt off 2.51pt}]  (277.67,141.33) -- (335,165.33) ;
\draw   (126.53,19.64) .. controls (125.65,16.95) and (128.52,14.29) .. (133.9,12.79) .. controls (139.29,11.28) and (146.25,11.2) .. (151.83,12.57) .. controls (153.81,11.01) and (157.43,9.93) .. (161.6,9.66) .. controls (165.77,9.39) and (169.99,9.96) .. (173,11.21) .. controls (174.68,9.79) and (177.99,8.84) .. (181.75,8.69) .. controls (185.51,8.54) and (189.18,9.21) .. (191.47,10.47) .. controls (194.52,8.97) and (199.36,8.34) .. (203.91,8.85) .. controls (208.45,9.36) and (211.89,10.92) .. (212.72,12.86) .. controls (216.45,13.29) and (219.56,14.37) .. (221.24,15.83) .. controls (222.92,17.29) and (223.01,18.99) .. (221.49,20.48) .. controls (225.16,22.49) and (226.02,25.16) .. (223.75,27.5) .. controls (221.47,29.84) and (216.4,31.5) .. (210.43,31.85) .. controls (210.39,34.05) and (207.52,36.07) .. (202.92,37.12) .. controls (198.33,38.18) and (192.72,38.12) .. (188.27,36.95) .. controls (186.38,39.58) and (181.05,41.52) .. (174.58,41.92) .. controls (168.11,42.33) and (161.67,41.13) .. (158.03,38.84) .. controls (153.58,39.97) and (148.23,40.29) .. (143.2,39.74) .. controls (138.17,39.19) and (133.88,37.81) .. (131.29,35.92) .. controls (126.74,36.14) and (122.34,35.15) .. (120.28,33.44) .. controls (118.21,31.73) and (118.92,29.66) .. (122.05,28.26) .. controls (117.99,27.26) and (115.92,25.27) .. (116.91,23.33) .. controls (117.91,21.4) and (121.75,19.95) .. (126.43,19.75) ; \draw   (122.05,28.26) .. controls (123.97,28.74) and (126.18,28.95) .. (128.4,28.88)(131.3,35.92) .. controls (132.25,35.87) and (133.18,35.77) .. (134.07,35.62)(158.03,38.84) .. controls (157.36,38.42) and (156.8,37.97) .. (156.36,37.5)(188.27,36.95) .. controls (188.62,36.47) and (188.84,35.98) .. (188.94,35.48)(210.43,31.85) .. controls (210.48,29.52) and (207.31,27.37) .. (202.29,26.35)(221.49,20.48) .. controls (220.67,21.28) and (219.43,21.98) .. (217.86,22.55)(212.72,12.86) .. controls (212.86,13.18) and (212.92,13.51) .. (212.91,13.83)(191.47,10.47) .. controls (190.71,10.85) and (190.09,11.26) .. (189.61,11.71)(173,11.21) .. controls (172.59,11.55) and (172.29,11.91) .. (172.1,12.28)(151.83,12.57) .. controls (153.01,12.86) and (154.11,13.21) .. (155.09,13.61)(126.53,19.64) .. controls (126.64,20.01) and (126.83,20.38) .. (127.09,20.74) ;
\draw    (152.67,40.33) -- (136.73,65.64) ;
\draw [shift={(135.67,67.33)}, rotate = 302.2] [color={rgb, 255:red, 0; green, 0; blue, 0 }  ][line width=0.75]    (10.93,-3.29) .. controls (6.95,-1.4) and (3.31,-0.3) .. (0,0) .. controls (3.31,0.3) and (6.95,1.4) .. (10.93,3.29)   ;
\draw    (185.33,39) -- (194.41,56.96) -- (200.1,68.22) ;
\draw [shift={(201,70)}, rotate = 243.19] [color={rgb, 255:red, 0; green, 0; blue, 0 }  ][line width=0.75]    (10.93,-3.29) .. controls (6.95,-1.4) and (3.31,-0.3) .. (0,0) .. controls (3.31,0.3) and (6.95,1.4) .. (10.93,3.29)   ;
\draw    (109.67,100.67) .. controls (107.57,110.47) and (101.67,108.67) .. (91.67,108.67) .. controls (81.67,108.67) and (66.18,115.42) .. (83.92,114.71) .. controls (100.16,114.06) and (89.82,125.17) .. (84.37,130.64) ;
\draw [shift={(83,132)}, rotate = 315] [color={rgb, 255:red, 0; green, 0; blue, 0 }  ][line width=0.75]    (10.93,-3.29) .. controls (6.95,-1.4) and (3.31,-0.3) .. (0,0) .. controls (3.31,0.3) and (6.95,1.4) .. (10.93,3.29)   ;
\draw    (423.67,191) .. controls (421.57,200.8) and (415.67,199) .. (405.67,199) .. controls (395.67,199) and (380.18,205.75) .. (397.92,205.04) .. controls (414.42,204.38) and (390.93,227.98) .. (383.81,235.04) ;
\draw [shift={(382.5,236.33)}, rotate = 315] [color={rgb, 255:red, 0; green, 0; blue, 0 }  ][line width=0.75]    (10.93,-3.29) .. controls (6.95,-1.4) and (3.31,-0.3) .. (0,0) .. controls (3.31,0.3) and (6.95,1.4) .. (10.93,3.29)   ;
\draw   (314.53,52.64) .. controls (313.65,49.95) and (316.52,47.29) .. (321.9,45.79) .. controls (327.29,44.28) and (334.25,44.2) .. (339.83,45.57) .. controls (341.81,44.01) and (345.43,42.93) .. (349.6,42.66) .. controls (353.77,42.39) and (357.99,42.96) .. (361,44.21) .. controls (362.68,42.79) and (365.99,41.84) .. (369.75,41.69) .. controls (373.51,41.54) and (377.18,42.21) .. (379.47,43.47) .. controls (382.52,41.97) and (387.36,41.34) .. (391.91,41.85) .. controls (396.45,42.36) and (399.89,43.92) .. (400.72,45.86) .. controls (404.45,46.29) and (407.56,47.37) .. (409.24,48.83) .. controls (410.92,50.29) and (411.01,51.99) .. (409.49,53.48) .. controls (413.16,55.49) and (414.02,58.16) .. (411.75,60.5) .. controls (409.47,62.84) and (404.4,64.5) .. (398.43,64.85) .. controls (398.39,67.05) and (395.52,69.07) .. (390.92,70.12) .. controls (386.33,71.18) and (380.72,71.12) .. (376.27,69.95) .. controls (374.38,72.58) and (369.05,74.52) .. (362.58,74.92) .. controls (356.11,75.33) and (349.67,74.13) .. (346.03,71.84) .. controls (341.58,72.97) and (336.23,73.29) .. (331.2,72.74) .. controls (326.17,72.19) and (321.88,70.81) .. (319.29,68.92) .. controls (314.74,69.14) and (310.34,68.15) .. (308.28,66.44) .. controls (306.21,64.73) and (306.92,62.66) .. (310.05,61.26) .. controls (305.99,60.26) and (303.92,58.27) .. (304.91,56.33) .. controls (305.91,54.4) and (309.75,52.95) .. (314.43,52.75) ; \draw   (310.05,61.26) .. controls (311.97,61.74) and (314.18,61.95) .. (316.4,61.88)(319.3,68.92) .. controls (320.25,68.87) and (321.18,68.77) .. (322.07,68.62)(346.03,71.84) .. controls (345.36,71.42) and (344.8,70.97) .. (344.36,70.5)(376.27,69.95) .. controls (376.62,69.47) and (376.84,68.98) .. (376.94,68.48)(398.43,64.85) .. controls (398.48,62.52) and (395.31,60.37) .. (390.29,59.35)(409.49,53.48) .. controls (408.67,54.28) and (407.43,54.98) .. (405.86,55.55)(400.72,45.86) .. controls (400.86,46.18) and (400.92,46.51) .. (400.91,46.83)(379.47,43.47) .. controls (378.71,43.85) and (378.09,44.26) .. (377.61,44.71)(361,44.21) .. controls (360.59,44.55) and (360.29,44.91) .. (360.1,45.28)(339.83,45.57) .. controls (341.01,45.86) and (342.11,46.21) .. (343.09,46.61)(314.53,52.64) .. controls (314.64,53.01) and (314.83,53.38) .. (315.09,53.74) ;
\draw    (340.67,73.33) -- (324.73,98.64) ;
\draw [shift={(323.67,100.33)}, rotate = 302.2] [color={rgb, 255:red, 0; green, 0; blue, 0 }  ][line width=0.75]    (10.93,-3.29) .. controls (6.95,-1.4) and (3.31,-0.3) .. (0,0) .. controls (3.31,0.3) and (6.95,1.4) .. (10.93,3.29)   ;
\draw    (373.33,72) -- (382.41,89.96) -- (388.1,101.22) ;
\draw [shift={(389,103)}, rotate = 243.19] [color={rgb, 255:red, 0; green, 0; blue, 0 }  ][line width=0.75]    (10.93,-3.29) .. controls (6.95,-1.4) and (3.31,-0.3) .. (0,0) .. controls (3.31,0.3) and (6.95,1.4) .. (10.93,3.29)   ;
\draw    (234.67,130.67) .. controls (232.57,140.47) and (226.67,138.67) .. (216.67,138.67) .. controls (206.67,138.67) and (191.18,145.42) .. (208.92,144.71) .. controls (225.16,144.06) and (214.82,155.17) .. (209.37,160.64) ;
\draw [shift={(208,162)}, rotate = 315] [color={rgb, 255:red, 0; green, 0; blue, 0 }  ][line width=0.75]    (10.93,-3.29) .. controls (6.95,-1.4) and (3.31,-0.3) .. (0,0) .. controls (3.31,0.3) and (6.95,1.4) .. (10.93,3.29)   ;
\draw    (116,101) .. controls (115.6,108.87) and (135.15,114.84) .. (129.11,119.84) .. controls (123.06,124.84) and (162.87,123.3) .. (141.93,133.15) .. controls (121,143) and (165.96,136.17) .. (153.14,146.13) .. controls (140.31,156.09) and (166,157) .. (163.08,148.3) .. controls (160.15,139.6) and (168.05,148.04) .. (177,146) .. controls (184.92,144.19) and (192.98,137.71) .. (197.45,134.2) ;
\draw [shift={(199,133)}, rotate = 143.13] [color={rgb, 255:red, 0; green, 0; blue, 0 }  ][line width=0.75]    (10.93,-3.29) .. controls (6.95,-1.4) and (3.31,-0.3) .. (0,0) .. controls (3.31,0.3) and (6.95,1.4) .. (10.93,3.29)   ;
\draw    (302,161) .. controls (301.6,168.87) and (321.15,174.84) .. (315.11,179.84) .. controls (309.06,184.84) and (348.87,183.3) .. (327.93,193.15) .. controls (307,203) and (351.96,196.17) .. (339.14,206.13) .. controls (326.31,216.09) and (352,217) .. (349.08,208.3) .. controls (346.15,199.6) and (354.05,208.04) .. (363,206) .. controls (370.92,204.19) and (378.98,197.71) .. (383.45,194.2) ;
\draw [shift={(385,193)}, rotate = 143.13] [color={rgb, 255:red, 0; green, 0; blue, 0 }  ][line width=0.75]    (10.93,-3.29) .. controls (6.95,-1.4) and (3.31,-0.3) .. (0,0) .. controls (3.31,0.3) and (6.95,1.4) .. (10.93,3.29)   ;

\draw (372.17,162.4) node [anchor=north west][inner sep=0.75pt]  [font=\small]  {$x_{1}$};
\draw (187.7, 100.07) node [anchor=north west][inner sep=0.7pt]  [font=\small]  {$x_{n-1}$};
\draw (103,75.4) node [anchor=north west][inner sep=0.75pt]  [font=\small]  {$x_{n}$};
\draw (342,172.73) node [anchor=north west][inner sep=0.75pt]  [font=\small]  {$u_{1}$};
\draw (144,106.73) node [anchor=north west][inner sep=0.75pt]  [font=\small]  {$u_{n-1}$};
\draw (168,49.73) node [anchor=north west][inner sep=0.75pt]  [font=\small]  {$\beta _{n}$};
\draw (244,76.4) node [anchor=north west][inner sep=0.75pt]  [font=\small]  {$\beta _{n-1}$};
\draw (362.67,120.4) node [anchor=north west][inner sep=0.75pt]  [font=\small]  {$\beta _{2}$};
\draw (66.67,132.4) node [anchor=north west][inner sep=0.75pt]  [font=\small]  {$\alpha _{n}$};
\draw (210,165.4) node [anchor=north west][inner sep=0.75pt]  [font=\small]  {$\alpha _{n-1}$};
\draw (391.67,231.73) node [anchor=north west][inner sep=0.75pt]  [font=\small]  {$\alpha _{1}$};
\draw (115.33,48.4) node [anchor=north west][inner sep=0.75pt]  [font=\small]  {$w$};
\draw (208,49.07) node [anchor=north west][inner sep=0.75pt]  [font=\small]  {$w$};
\draw (303.33,81.4) node [anchor=north west][inner sep=0.75pt]  [font=\small]  {$w$};
\draw (398.33,96.4) node [anchor=north west][inner sep=0.75pt]  [font=\small]  {$w$};
\draw (320.33,209.4) node [anchor=north west][inner sep=0.75pt]  [font=\small]  {$v_{1}$};
\draw (132.33,153.4) node [anchor=north west][inner sep=0.75pt]  [font=\small]  {$v_{n-1}$};

\end{tikzpicture}

\caption{Diagram of the water flow network in Section~\ref{sect_examples}.}
    \label{figure_rain}
\end{figure}
The system dynamics corresponding to Figure~\ref{figure_rain} are modeled as follows
\begin{align}\label{dyn_line_shape}
    &\dot{x}_1=-(\alpha_1)x_1+\beta_{2}x_{2}-u_1+ v_1, \notag\\
    &\dot{x}_i=-(\alpha_i+\beta_i)x_i+\beta_{i+1}x_{i+1}-u_i+u_{i-1} + v_i-v_{i-1}  \notag\\
    &\dot{x}_n=-(\alpha_n+\beta_n)x_n+u_{n-1} - v_{n-1}
\end{align}
$i=2, \dots, n-1$, where $x_i$ represents the water volume at each section $i$. Assume that the network is homogeneous, i.e. $\alpha_i=\alpha$ for all $i=1,\dots, n$,  $\beta_i=\beta$ for all $i=2,\dots, n$. Then the system matrices $A\in \mathbb R^{n \times n}$, $B\in \mathbb R^{n\times m}$ with $m=n-1$ and initial state $x(0)=x_0$ are 
\begin{align}\label{A_B_H_line}
&\scalebox{0.8}{$A=\begin{bmatrix}
-\alpha &\beta   & \cdots &0 \\ 
0 & -(\alpha+\beta)   & \cdots &0 \\  
\vdots & \vdots & \ddots & \beta \\
0 & 0 & 0 &-(\alpha+\beta)
\end{bmatrix}; \hspace{1mm} B=\begin{bmatrix}
-1 & 0 & \dots & 0 \\
1 & -1 & \dots & 0 \\
0 & \ddots & \ddots & \vdots \\
 \vdots & \ddots & \ddots & -1 \\
0 & 0 \hspace{2mm}\cdots& 0 & 1
\end{bmatrix}$;} \notag \\
&\scalebox{0.8}{$H=-B; \hspace{47mm} x_0=\mathbb{1}$.}
\end{align}
It is worth noting that the model’s formulation enables scalability, allowing for the analysis of extensive river networks and accommodating changes in topology, such as the addition of control stations or the influence of evolving environmental disturbances.

\subsection{Leakages: Load Disturbances}
Initially, assume there is no additional rain input, i.e. $w=0$. The control and disturbance actions in this network system are constrained by $\left|u \right|\leq Ex$, $\left|v \right|\leq Gx$ with $E,G \geq 0$. Here, $G$ represents the slope of the water flow network, where a higher altitude corresponds to a greater leakage disturbance. The choice of parameters $E, s, r, \delta$ are designed to satisfy the assumptions of Theorem~\ref{theorem_finite}, see~\eqref{ass1_ex},~\eqref{ass2_ex}:
\begin{align}\label{E_G_s_r_delta}
&\scalebox{0.8}{$E=  \begin{bmatrix}
0& \zeta_u & 0 & \cdots &  0 \\
0 & 0& \zeta_u & \cdots  & 0 \\
\vdots & \vdots & \ddots & \ddots & \vdots \\
0 & 0 & \dots  &  0& \zeta_u
\end{bmatrix}; \hspace{2mm}  G = \frac{1}{n} \begin{pmatrix}
\zeta_v & 0 &  \dots &  0 & 0 \\
0 & 2\zeta_v &  \dots &  0 & 0 \\
\vdots & \vdots & \ddots & \vdots & \vdots  \\
0 & 0 &  \dots &  (n-1)\zeta_v & 0 
\end{pmatrix};$} \notag\\ 
&\scalebox{0.9}{$s= \rho_s \begin{bmatrix}
1 \hspace{1mm} 0 \hspace{1mm} \dots \hspace{1mm} 0 \hspace{1mm} 0\\
\end{bmatrix}; \hspace{5mm} r= \rho_u \begin{bmatrix}
2/n \hspace{1mm} \dots \hspace{1mm} (n-1)/n \hspace{1mm} 0 
\end{bmatrix}^{\top};$}
\end{align}
and $\delta=\rho_v\mathbb{1}$. Here $\zeta_u$ and $\zeta_v$ scale, respectively, the control and disturbance capacity. $\rho_s$, $\rho_u$, $\rho_v$ are scaling parameters for the penalties on state, control, and disturbance terms in the cost function. The sparsity structure of $E$ relies on the Metzler condition, ensuring a viable setup under the given constraints.

The optimal control problem setup~\eqref{minimaxProb_finite} of this example is
\begin{align}\label{minimax_prob_setup_inf_example}
    &\underset{\mu}{\min} \hspace{1mm} \underset{v}{\sup} \hspace{2mm} \int_{0}^{T} \left [ s^{\top}x(\tau)+ r^{\top}u(\tau)- \delta^{\top}v(\tau)\right ]d\tau \notag \\
            &\mathrm{Subject \hspace{2mm} to} \notag \\
        &~~~~~~~~\dot{x}(t)=Ax(t)+Bu(t)+Hv(t) \\
        &~~~~~~~~x(0)=x_{0}, \hspace{2mm} u(t)=\mu(x(t)), \notag\\
        &~~~~~~~~\left | u \right | \leq Ex, \hspace{2mm} \left | v \right | \leq G x. \notag
\end{align}
In this example, the objective is to optimize the performance at the lowest point of the river. The objective function is defined as the integral of the financial cost $z$ over time. This cost includes the state penalty vector $s$ representing the cost associated with maintaining $x(t)\geq 0$, $r$ accounting for the operation of the dams over time, and $\delta$ reflecting the savings due to natural water dissipation.

The choice of $\zeta_u$, $\zeta_v$, $\rho_s$, $\rho_u$, and $\rho_v$ depends on conditions~\eqref{ass_a},~\eqref{ass_s}. Note that, under the assumption that the disturbance input gradient $-\delta+H^{\top}p(t)$ does not change sign for all $t$, it is possible that the optimal disturbance does not violate the Metzler condition of the closed-loop system matrix. In this case, it suffices that $A-|B|E$ is Metzler, given by
\begin{align}\label{ass1_ex}
    \scalebox{0.9}{$A-\left| B\right|E 
=\begin{bmatrix}
-\alpha & \beta - \zeta_u   & \cdots &  0\\
0 & -(\alpha + \beta) -\zeta_u &\cdots &  0\\
\vdots & \vdots & \ddots & \vdots \\
0 & 0 & \cdots &  \beta- \zeta_u  \\
0 & 0 & \cdots & -(\alpha + \beta)- \zeta_u 
\end{bmatrix}$}
\end{align}
Under this relaxed condition, it is necessary that $\beta \geq\zeta_u$.
Moreover, it is also required
\begin{align}\label{ass2_ex}
    \scalebox{0.9}{$s -\left( E ^{\top} \left | r \right | - G^{\top} \left | \delta \right | \right)
=\begin{bmatrix}
\rho_s + \frac{1}{n}\zeta_v \rho_v\\
0- \frac{2}{n}\left(\zeta_u \rho_u - \zeta_v \rho_v\right)\\
\vdots \\
0- \frac{n-1}{n}\left(\zeta_u \rho_u - \zeta_v \rho_v\right)\\
0
\end{bmatrix}\geq 0.$}
\end{align}
Hence, for all $\rho_s\geq 0$ it is sufficient that, $ \zeta_u \rho_u - \zeta_v \rho_v \leq 0.$
 
Set $\alpha=3$, $\beta=10$, $\zeta_u=6$, $\zeta_v=3$, $\rho_s=1$, $\rho_u=0.17$, $\rho_v=0.35$ and assume the river flow has $n= 100$ sections. Note that, under these parameters, the state matrix has a Perron-Frobenius eigenvalue of $\lambda_p=-3$.
Since $|v|\le Gx$, the contribution $Hv$ of the disturbance to the dynamics of $x$ is upper-bounded by $|H|Gx$. The dynamics $\dot x = (A+|H|G)x + Bu$ therefore represent an overestimation of the effect of the disturbance. Next, it is demonstrated that the overestimated disturbance setting can also be attenuated by this worst-case disturbance model.
\subsubsection{Detectability} Let us verify the detectability as presented in Section~\ref{Stab_Detect_subsec} for the overestimated disturbance setting, with $\tilde{A}= A+|H|G$ and $\tilde{s}= s^{\top} + |\delta^{\top}|G$. The detectability of the pair $(\tilde{s}-r^{\top}K, \tilde{A}-BK)$ is established via  Theorem~\ref{thm_detect_cond}. However, the sufficient condition provided in Corollary~\ref{cor_suff_cond_detec}  is overly restrictive for this example, as \eqref{ass2_ex} is not strictly positive. Since $\lambda_p \approx 0.19$ it suffices to choose $h=20$ to verify the inequality
 \begin{align*}
         \frac{1}{h}(\tilde{s}^{\top}-|r^{\top}|E) \sum_{i=0}^{n-1}\left(\frac{1}{h}(\tilde{A}-|B|E)+I \right)^i>0.
     \end{align*}
Therefore, the pair $(\tilde{s}-r^{\top}K, \tilde{A}-BK)$ is detectable for any $K$ satisfying $|K| \leq E.$
\subsubsection{Simulations}
\begin{figure}[h]
    \centering
    \includegraphics[scale=0.5]{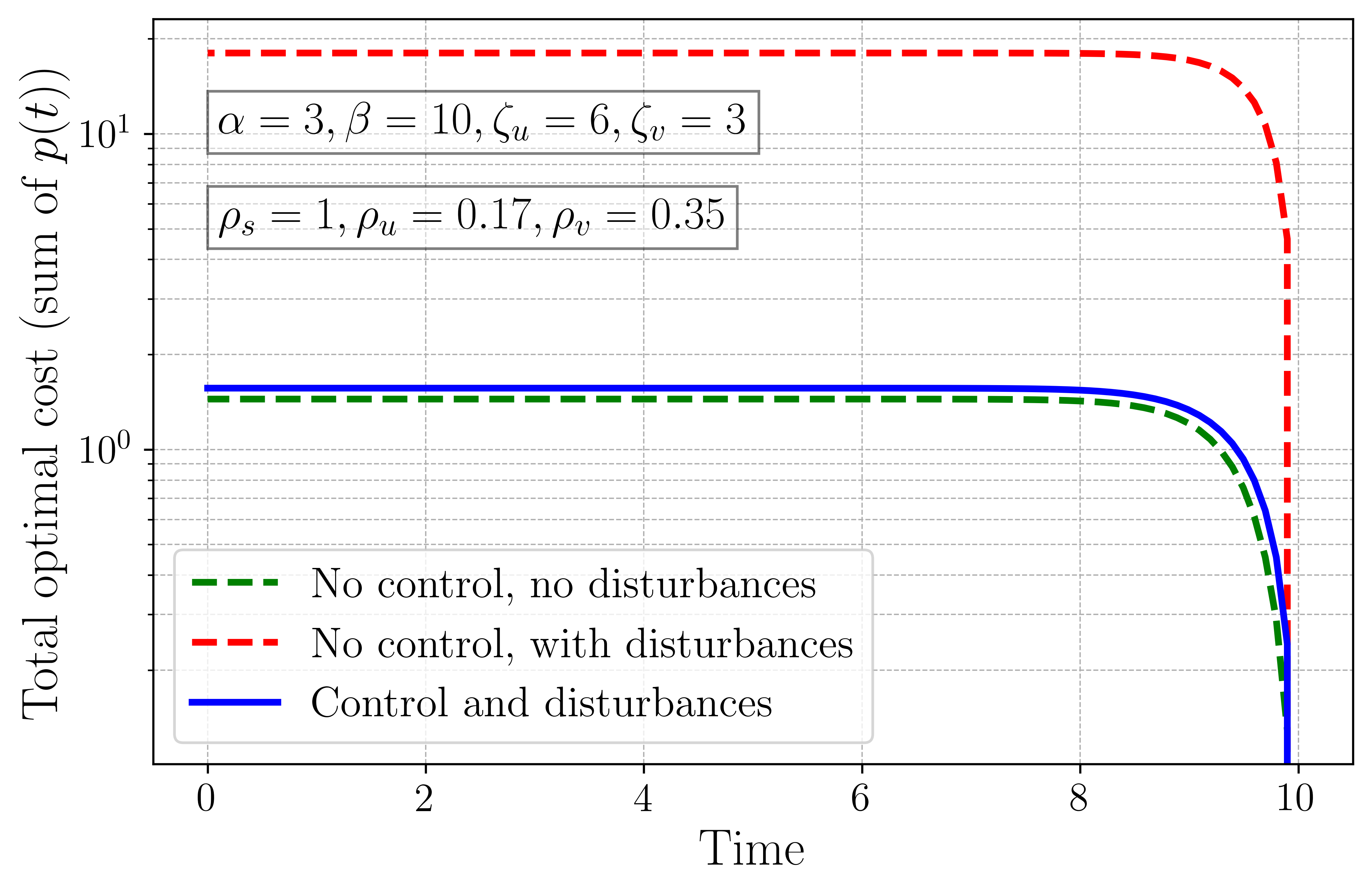}
    \caption{Evolution of the optimal cost $p(t)^{\top}\mathbb{1}$ over the time interval $t \in \left[0, 10\right]$ plotted on a logarithmic scale.}
    \label{finitefig1}
\end{figure}

In Figure~\ref{finitefig1}, the optimal cost evolution shows that, despite the optimal disturbances causing an increase in the cost, the controller derived from our minimax framework successfully reduces the cost to the level observed in the absence of disturbances. This demonstrates that, under pertinent assumptions and proper tuning of the different optimization variables, the controller’s effectiveness allows the system to maintain performance comparable to the disturbance-free scenario. 
It is noted that $-\delta + H\T p(t) < 0$ for all $t$, and therefore the worst-case disturbance is given by $v^* = -Gx$. Closing the disturbance loop gives the dynamics $\dot x = (A-HG)x + Bu$, which is Metzler by \eqref{ass1_ex} and does not destabilize the system dynamics.

\begin{figure}[h]
    \centering
    \begin{minipage}{0.75\textwidth}
        \centering
        \includegraphics[width=\textwidth]{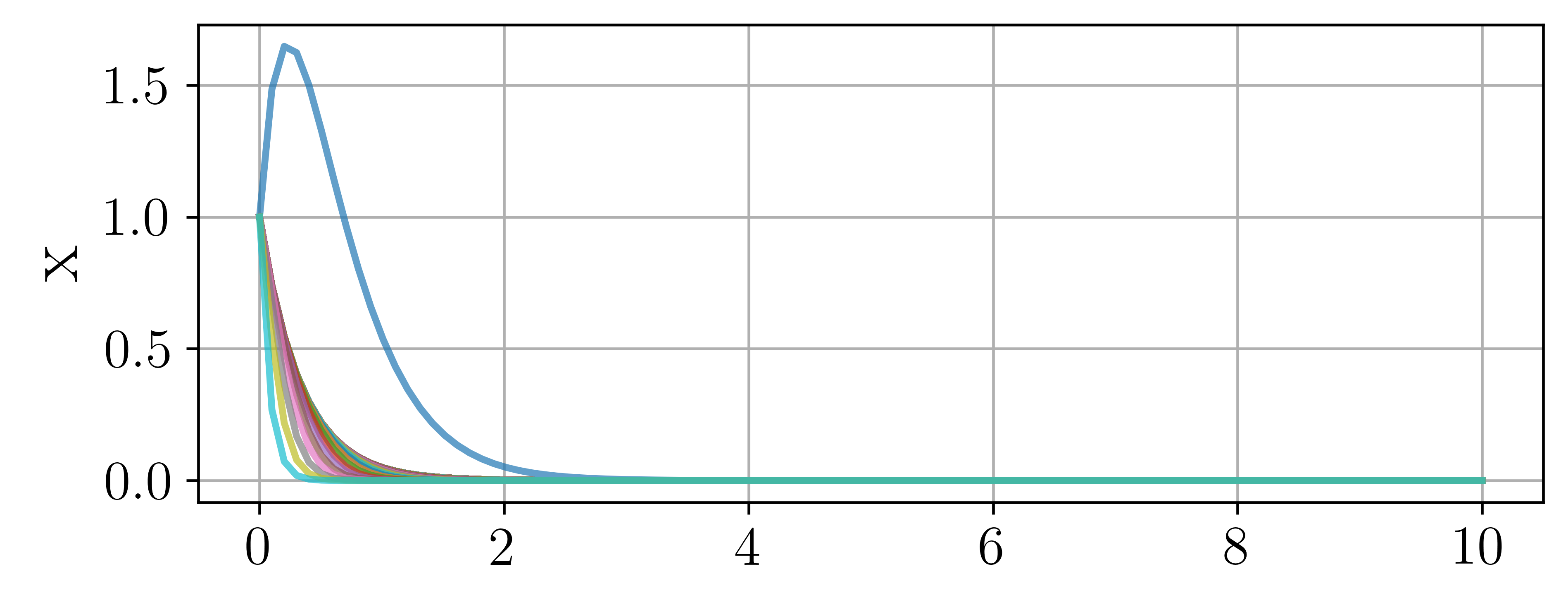}\\
        \small (a)
    \end{minipage}
    \begin{minipage}{0.75\textwidth}
        \centering
        \includegraphics[width=\textwidth]{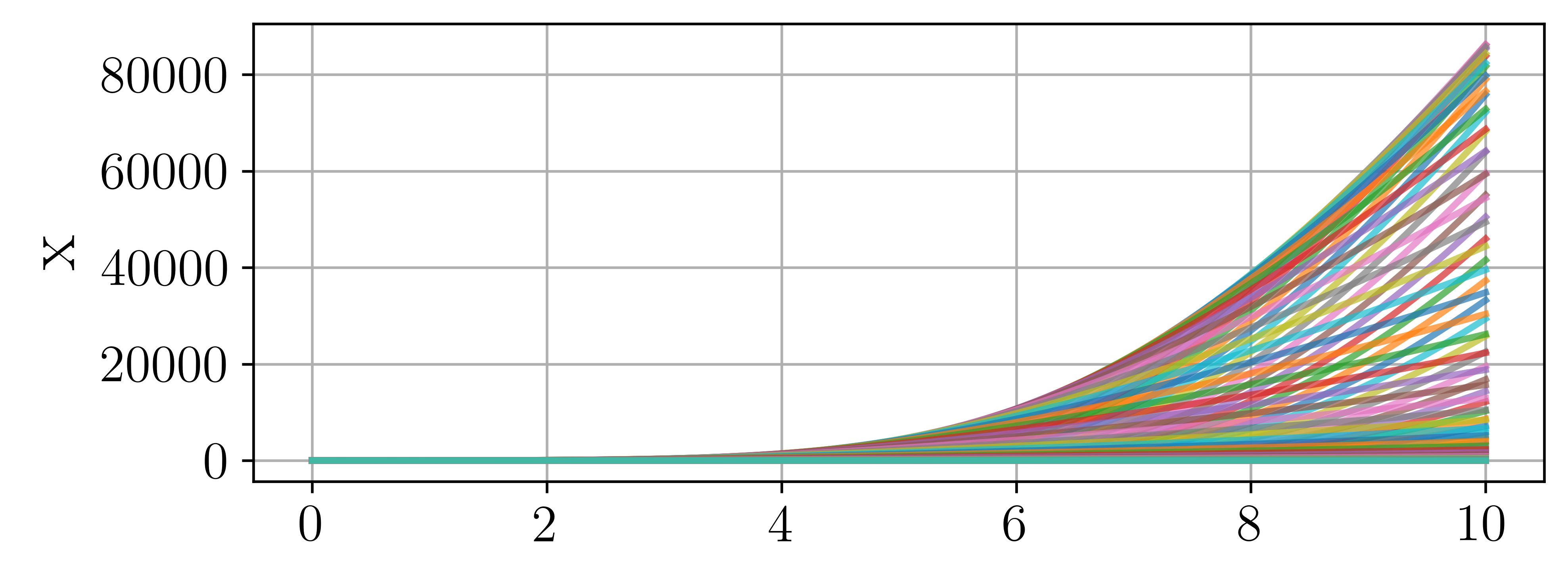}\\
        \small (b)
    \end{minipage}
    \begin{minipage}{0.75\textwidth}
        \centering
        \includegraphics[width=\textwidth]{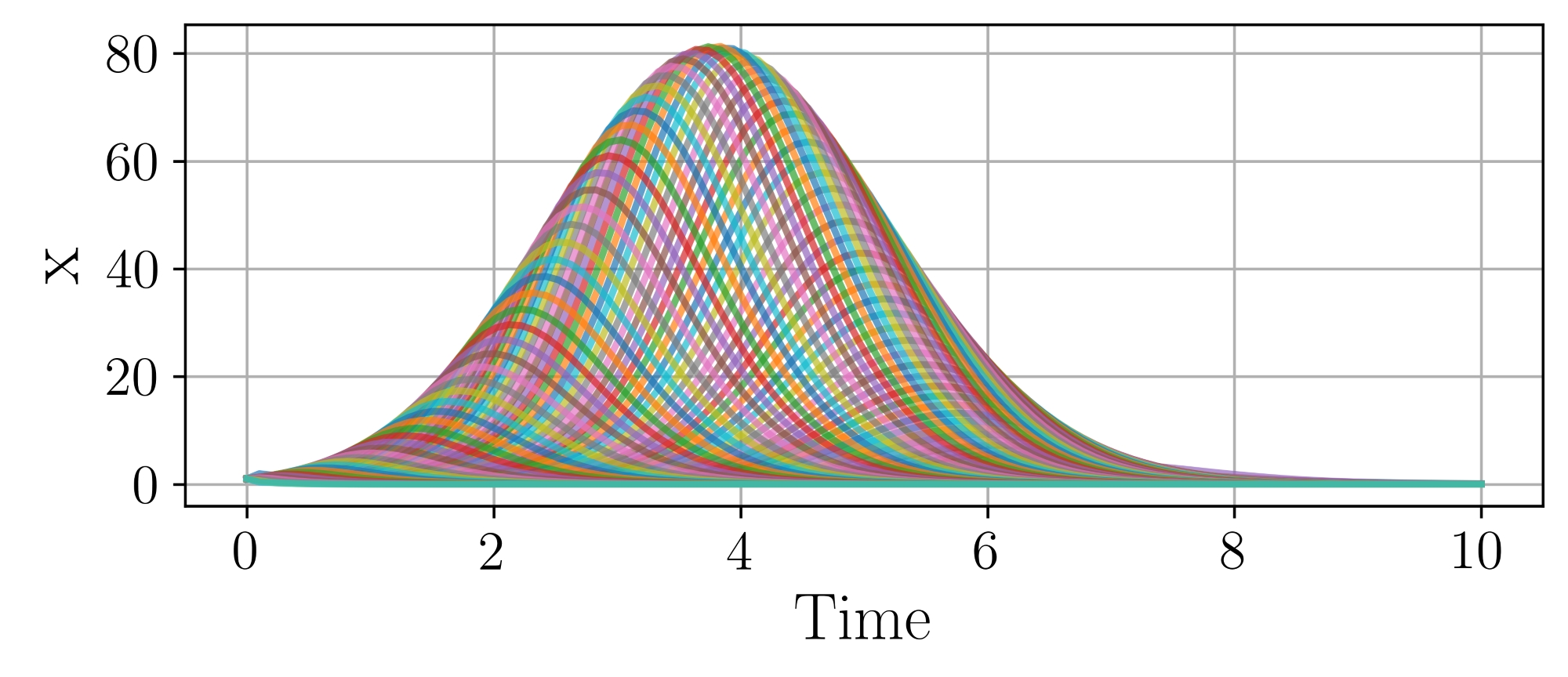}\\
        \small (c)
    \end{minipage}
    \caption{State trajectories $x(t)=e^{\tilde{A}t}x_0$ over the time horizon $t \in \left[0, 10\right]$, under different configurations. (a) Open-loop evolution with $\tilde{A}=A$; (b) open-loop dynamics in the presence of disturbances, $\tilde{A}=A+\left|H \right|G$; and (c) closed-loop behavior in the presence of both disturbances and control, given by $\tilde{A}=A+\left|H \right|G-BK$.}
    \label{finitefig2}
\end{figure}

Conversely, Figure~\ref{finitefig2} demonstrates that the controller not only counteracts the increase of the optimal cost introduced by the optimal disturbance policy, and also compensates for the overestimated disturbance contribution in the system dynamics, given by $\tilde A=A+|H|G$. This contribution, while not realized by the optimizing disturbance, represents an upper-bound direction that drives the open-loop system unstable, with $\lambda_p\approx 0.19$. The figure illustrates that the same minimax controller is able to stabilize the system even in this scenario, which demonstrates its robustness. %

\begin{figure}[h]
    \centering
    \includegraphics[scale=0.5]{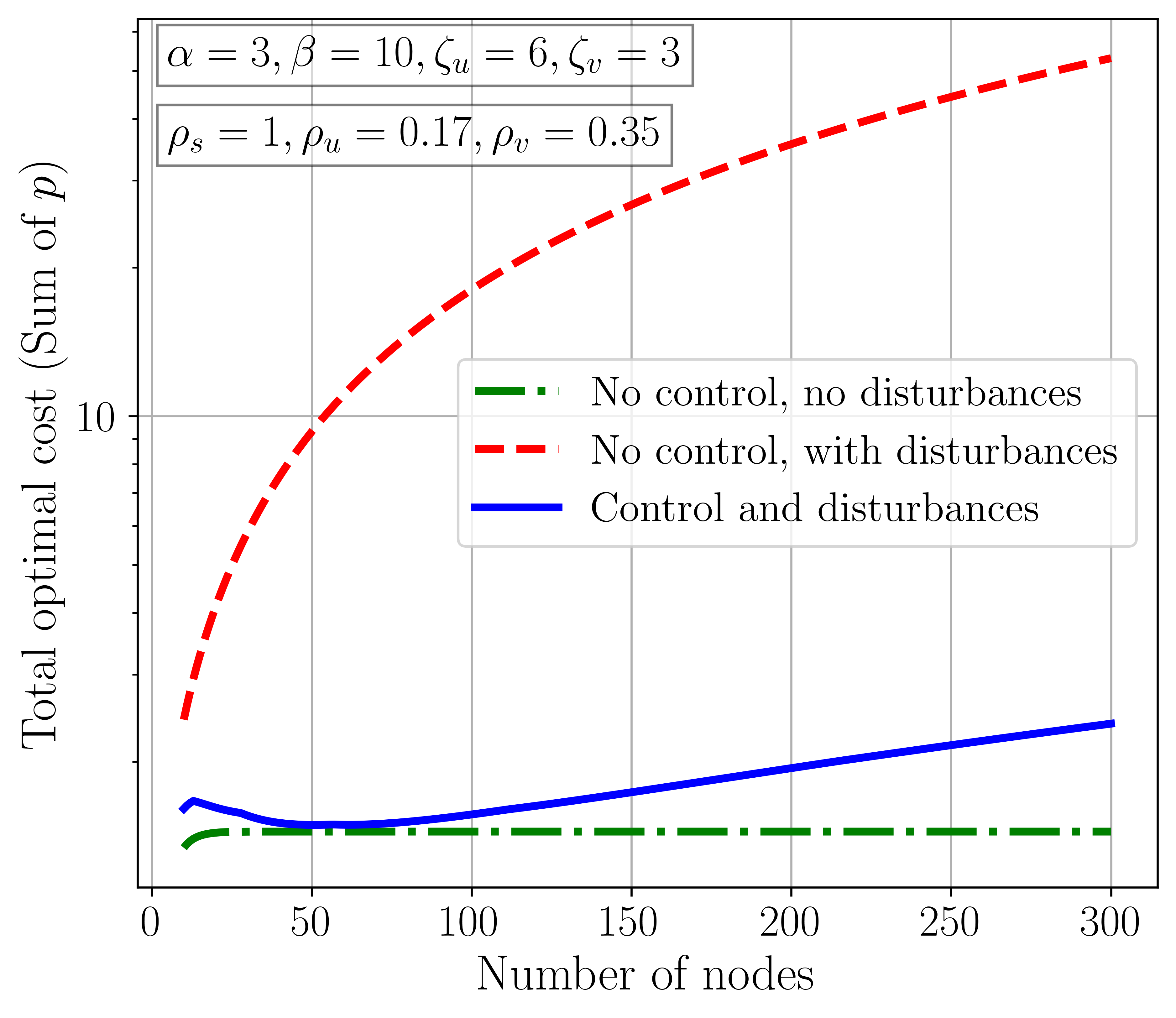}
    \caption{Optimal cost evolution $p^{\top}\mathbb{1}$ when the time horizon tends to infinity, as the water-flow network is scaled from $n=2$ to $200$ segments on a logarithmic scale.}
    \label{inffig}
\end{figure}

Finally, Figure~\ref{inffig} illustrates the behavior of the total optimal cost when $T \rightarrow \infty$ and as the network size increases. As the number of nodes grows, the cost in the presence of disturbances rises significantly, indicating how disturbances escalate the cost in larger networks. In contrast, when the controller derived from the minimax framework is introduced, it significantly reduces again the cost.
\subsection{Rain: Positive Unconstrained Disturbance.}
Assume now that it starts raining. Let the rain be interpreted as a worst-case positive, unconstrained disturbance $w\geq 0$ affecting all nodes homogeneously, i.e. $F=\mathbb{1}$ and $\gamma$ be the permeability of the soil. Suppose this disturbance persists over a time horizon of $T=24$. An upper bound on the amount of rain that the feedback controller, derived from the minimax framework, can compensate for is studied. To establish this bound, condition~\eqref{gamma_ass_finite} from Theorem~\ref{theorem_finite} is applied. This condition provides a measure of the system’s capacity to handle the disturbance, with $\boldsymbol{\gamma} \geq F^{\top}p(0)$ representing the minimum $L_1$-induced gain required to counteract the worst-case rain disturbance and maintain system stability.
The extended optimal control problem becomes
\begin{align}\label{minimax_prob_setup_inf_example2}
    &\underset{\mu}{\min} \hspace{1mm} \underset{v}{\sup} \hspace{2mm} \int_{0}^{24} \left [ s^{\top}x(\tau)+ r^{\top}u(\tau)- \delta^{\top}v(\tau)-\boldsymbol{\gamma}^{\top}w\right ]d\tau \notag \\
            &\mathrm{Subject \hspace{2mm} to} \notag \\
        &~~~~~~~~\dot{x}(t)=Ax(t)+Bu(t)+Hv(t)+\mathbb{1}w \\
        &~~~~~~~~x(0)=x_{0}, \hspace{2mm} u(t)=\mu(x(t)), \notag\\
        &~~~~~~~~\left | u \right | \leq Ex, \hspace{2mm} \left | v \right | \leq G x, \hspace{2mm} w\geq 0. \notag
\end{align}
Assuming the same parameter choice as before, it is obtained that $F^{\top}p(0) \approx 1.56$. Therefore, if the optimal control problem~\eqref{minimax_prob_setup_inf_example} has a finite solution $p^{\top}(0)x_0$ and $\gamma \geq 1.56$ then the extended optimal control problem~\eqref{minimax_prob_setup_inf_example2} also has a finite solution, and its solution coincides with that of~\eqref{minimax_prob_setup_inf_example}. This result establishes that the controller derived for the original problem is robust enough to handle the additional rain disturbance, provided that the gain condition is satisfied.

\section{CONCLUSIONS}\label{sect7_conclusions}

This work provides an extensive framework for minimax problems for positive systems that considers bounded and unconstrained worst-case disturbances.  
This work significantly extends prior results in discrete time~\cite{AlbaEmmaAnders} by deriving explicit solutions for the continuous-time setting. Notably, despite minimal assumptions on the controller, the resulting optimal policy is linear, with its sparsity structure inherently determined by the problem constraints.
Moreover, the \nobreak{Linear Regulator} problem and its linear programming formulation were introduced and studied. Stabilizability and detectability aspects of the system dynamics are studied and related to the existence and equivalence of solutions of the HJI equation and the linear program. Additionally, sufficient detectability conditions for the uniqueness of a stabilizing policy when the HJI equation has a finite solution are derived.
It was shown that the linear program always finds a stabilizing policy, if it exists.
In addition, the $L_1$-induced gain of the system with respect to unbounded disturbance was formulated, and a tight bound for the finiteness of the cost under this disturbance was presented.
The approach was motivated by an example that considers a large-scale water management network.

The presented theoretical framework enables efficient scaling to large dynamical systems. This efficiency stems partly from the sparsity of the optimal policy, which is a direct result of the constraints imposed by the problem statement and the optimization process. This aspect is currently a subject of ongoing research. Further analysis of network structural properties that enable particular efficiency and scalability is subject of ongoing work. 
Another area of future work involves improving our fixed-point iteration for computing the solution of the algebraic HJI equation, since a linear program has not yet been formulated.
\appendix\label{Appendix}
In this Appendix, we build upon existing results in continuous-time optimal control ~\cite{Subbotin},~\cite{Bertsekas_dynprograndoptcontrol, basar}, which serve as an auxiliary foundation for proving the main theorems of this paper.

Consider the continuous-time dynamics
\begin{align}\label{11.1}
    \dot{x}(t)=f(t,x(t), u(t), v(t)), \hspace{2mm} 0\leq t\leq T
\end{align}
where $T\in \mathbb R_+$ is the terminal time horizon, $x(t)$ represents the vector of $n-$dimensional state variables at time $t$, $u(t)\in U$ the $m-$dimensional control variable vector at time $t$, $v(t)\in \Omega$ the $c-$dimensional disturbance at time $t$. The sets $U$, $\Omega$ are, respectively, the control and disturbance constrained sets and are both compact. Both the control and disturbance policies are chosen depending on the current position $(t, x(t))$. 
Consider a differential game governed by the differential equation~\eqref{11.1}, where the control player tries to minimize the cost 
\begin{align}\label{11.2}
    J(t_0,x(\cdot), u(\cdot), v(\cdot)):=\phi(x(T))+\int_{t_0}^{T} g(x(\tau), u(\tau), v(\tau)) \hspace{1mm} d \tau
\end{align}
where $t_0\in [0, T]$ is an initial instant of time. On the contrary, the disturbance player wishes to maximize this functional. 
The functions $f,g,\phi$ satisfy the following assumptions
\begin{ass}\label{A1}
    The functions $f:[0,T]\times \mathbb R^n \times U\times \Omega \mapsto \mathbb R^n$, $g: [0, T]\times \mathbb R^n \times U\times \Omega \mapsto \mathbb R$, and $\phi:\mathbb R^n \mapsto \mathbb R$ are continuous and satisfy
    \begin{align}\label{11.3}
        \left\|f(t,x,u,v)  \right\|&\leq (1+\left\|x  \right\|)\kappa_f \notag\\
        \left|g(t,x,u,v)  \right|&\leq (1+\left\|x  \right\|)\kappa_g\\
        |\phi(x)|&\leq (1+\left\|x  \right\|)\kappa_{\phi}\notag
    \end{align}
    for all $(t,x,u,v)\in [0, T]\times \mathbb R^n\times U\times \Omega$, where $\kappa_f$, $\kappa_g$, $\kappa_{\phi}$ are positive scalars.
\end{ass}
\begin{ass}\label{A2}
The functions $f$ and $g$ satisfy the Lipschitz condition in the variable $x$
\begin{align}\label{11.4}
    &\left\|f(t, x+y, u, v)-f(t,x,u,v)  \right\| \\
    &~~+\left| g(t,x+y,u,v)-g(t,x,u,v) \right|\leq \lambda \left\| y \right\|\notag
\end{align}
for all $(t,x,u,v)\in [0,T]\times \mathbb R^n\times U\times \Omega$, $y\in \mathbb R^n$.
\end{ass}
\begin{ass}\label{A3}
For any $q\in \mathbb R^n$ and $(t,x)\in [0, T]\times \mathbb R^n$ the equality 
\begin{align}
    &\min_{u\in U} \hspace{1mm}\max_{v \in \Omega} \left\{ q^{\top}f(t,x,u,v)-g(t,x,u,v)\right\}\\
    &~~~~~=  \max_{v \in \Omega} \hspace{1mm} \min_{u\in U} \left\{ q^{\top}f(t,x,u,v)-g(t,x,u,v) \right\}\notag
\end{align}
is valid.
\end{ass}

In this context, it is possible to formulate the differential game as the general, continuous-time, minimax optimal control problem with continuous constraints and nonnegative final conditions $\phi(x(T))$ as 
\begin{align}\label{gral_minimax_prob_setup_cont_finite}
    \underset{\mu}{\mathrm{inf}}\hspace{1mm}\underset{\omega}{\mathrm{\max}}& \hspace{3mm}\left[\phi(x(T))+\int_{0}^{T} g(x(\tau), u(\tau), v(\tau)) \hspace{1mm} d \tau \right] \notag\\
    \mathrm{subject}& \hspace{1mm} \mathrm{to}  \\
    &\dot{x} = f(x(t), u(t), v(t)) \notag \\
    &  x(t)\in X \subset \mathbb{R}^{n}; \hspace{2mm} x(0)=x_{0}\hspace{1mm}; \hspace{2mm} u(t)=\mu(x(t)) \notag\\
    &u(t)\in U(x(t)) \subset \mathbb{R}^{m}; \hspace{2mm} v(t)\in \Omega(x(t))\subset \mathbb{R}^{c}\notag
\end{align}
\begin{lem}\label{LEMA_cont_minimax_finite}
    Under Assumptions~\ref{A1},~\ref{A2},~\ref{A3} there exists the optimal value of the differential game~\eqref{gral_minimax_prob_setup_cont_finite}. The value function coincides with the minimax solution of the Cauchy problem 
    \begin{align}\label{HJB_minimax_cont_finite}
        -\nabla_t J(t,x)&=\min_u \hspace{1mm}\max_v \left[g(x,u,v)+\nabla_x J(t,x)^{\top}f(x,u,v) \right] \notag\\
        J(T,x)&=\phi(x(T)).
    \end{align}
\end{lem}
\begin{proof}
    The proof of this theorem can be found in~\cite[Sec. 12]{Subbotin}.
\end{proof}

 \bibliographystyle{abbrv} 
 \bibliography{cas-refs}
\end{document}